\documentclass[11pt, leqno, a4paper, twoside]{amsart}

\usepackage[utf8]{inputenc}
\usepackage{lineno}
\usepackage[numbers,sort&compress]{natbib}
\usepackage[draft]{changes}
\definechangesauthor[name=Paola, color=magenta]{PA}
\usepackage{amsmath}
\usepackage{siunitx}
\usepackage{amsfonts,bm}
\usepackage{dsfont}
\usepackage[colorlinks=true]{hyperref}
\usepackage{xcolor}
\usepackage{stmaryrd}
\usepackage{algorithm}
\usepackage{amssymb}
\usepackage{verbatim}
\usepackage{mathrsfs}
\usepackage{mathtools,amsmath}
\usepackage[margin=1.15in]{geometry}
\usepackage{tikz}
\usepackage{pgfplots}
\usepackage{capt-of}
\usepackage{mdframed}
\usepackage{multirow}
\usepackage{colortbl}
\numberwithin{equation}{section}

\usepackage[noend]{algpseudocode}
\usepackage[capitalise, noabbrev]{cleveref}
\usepackage{caption}
\usepackage{url}
\newtheorem{remark}{Remark}

\newcommand*{\llbrace}{\{\mskip-5mu\{}
\newcommand*{\rrbrace}{\}\mskip-5mu\}}
\renewcommand{\vec}{\underline}
\newcommand{\vecc}[1]{\mathbf{\underline{#1}}}

\newcommand{\mat}[1]{\textbf{#1}}
\newcommand{\tens}[1]{\boldsymbol{#1}}
\renewcommand{\epsilon}{\varepsilon}

\newcommand{\jump}[1]{\left[\mskip-2mu\left[#1\right]\mskip-2mu\right]}
\newcommand{\mean}[1]{\llbrace#1\rrbrace}

\newcommand{\pgv}{\textup{PGV}_{\textup{gmh}}}
\newcommand{\pgu}{\textup{PGU}_{\textup{gmh}}}

\definecolor{gray}{rgb}{0.5 0.5 0.5}
\definecolor{sgreen}{HTML}{40A441}
\definecolor{syellow}{HTML}{DFDE62}
\definecolor{spink}{HTML}{6D2A6F}
\definecolor{sred}{HTML}{B84041}
\definecolor{sblue}{HTML}{4040B2}
\definecolor{sgray}{HTML}{B6B6B6}
\definecolor{lltgray}{rgb}{0.8 0.8 0.8}
\definecolor{ltgreen}{rgb}{0.8,1,0.7}
\definecolor{dkgreen}{rgb}{0.4,0.6,0.35}
\definecolor{MATgreen}{HTML}{4CAF50}
\definecolor{MATblue}{HTML}{2196F3}

\newcommand{\kommentar}[1]{}

\title[Simulation of seismic response of dams]{Elasto-acoustic modelling and simulation for the seismic response of structures: The case of the Tahtal{\i} dam in the 2020 $\dot{\textup{I}}$zmir earthquake}


\author[I. Mazzieri, M. Muhr, M. Stupazzini and B. Wohlmuth]{Ilario Mazzieri$^1$, Markus Muhr$^{\ast,2}$, Marco Stupazzini$^{3}$, and Barbara Wohlmuth$^2$}

\email{\href{ilario.mazzieri@polimi.it}{ilario.mazzieri@polimi.it}}
\email{\href{mailto:muhr@ma.tum.de}{muhr@ma.tum.de}}
\email{\href{mailto:MStupazzini@munichre.com}{MStupazzini@munichre.com}}
\email{\href{mailto:wohlmuth@ma.tum.de}{wohlmuth@ma.tum.de}}

\thanks{$^*$Corresponding author: Markus Muhr, \href{mailto:muhr@ma.tum.de}{muhr@ma.tum.de}}

\begin{document}
\maketitle
\vspace*{-4mm}
\begin{center}
{\footnotesize
   $^1$MOX,  Dipartimento di Matematica, Politecnico di Milano, Milano, Italy  \\
   $^2$Department of Mathematics, Technical University of Munich, Garching, Germany \\
   $^3$Munich RE, M\"unchener R\"uckversicherungs-Gesellschaft, Munich, Germany
}
\end{center}

\vspace{8mm}


\begin{abstract}
As a mean to assess the risk dam structures are exposed to during earthquakes, we employ an abstract mathematical, three dimensional, elasto-acoustic coupled wave-propagation model taking into account (i) the dam structure itself, embedded into (ii) its surrounding topography, (iii) different material soil layers, (iv) the seismic source as well as (v) the reservoir lake filled with water treated as an acoustic medium. As a case study for extensive numerical simulations we consider the magnitude 7 seismic event of the 30$^{\rm th}$ of October 2020 taking place in the Icarian Sea (Greece) and the Tahtal{\i} dam around 30\,km from there (Turkey). A challenging task is to resolve the multiple length scales that are present due to the huge differences in size between the dam building structure and the area of interest, considered for the propagation of the earthquake. Interfaces between structures and highly non-conforming meshes on different scales are resolved by means of a discontinuous Galerkin approach. The seismic source is modeled using inversion data about the real fault plane. Ultimately, we perform a real data driven, multi-scale, full source-to-site, physics based simulation based on the discontinuous Galerkin spectral element method, which allows to precisely validate the ground motion experienced along the Tahtal{\i} dam, comparing the synthetic seismograms against actually observed ones. A comparison with a more classical computational method, using a plane wave with data from a deconvolved seismogram reading as an input, is discussed.
\end{abstract}



\section{Introduction}
With the continuous growth of computational power in the last decades, physics based simulation (PBS) emerged as an aspiring, alternative approach to ground motion prediction equations (GMPEs), which has already been applied to seismic scenarios at various sites including the United States \cite{PBS_Memphis, SCECCyberShake}, Japan \cite{PBS_Japan, Iwaki2016ValidationOT}, New-Zealand \cite{Bradley2020, Bradley2021}, Turkey \cite{IstanbulCaseStudy}, China \cite{MOXRiskAssessment}, the Netherlands \cite{Groeningen}, Italy \cite{ItalyCaseStudy}. PBS aims at describing, as reliably as possible, the seismic wave propagation problem and therefore it is crucial, on the one hand, to properly characterize the mechanical properties of the different portion of the computational domain and, on the other, to have a reliable seismic excitation source, see, e.g. \cite{Petersson2018, Peltieshpc, Galvez2014,Zelst2019, IstanbulCaseStudy, Groeningen}. Because of the intrinsically high epistemic uncertainties involved in the construction of 3D numerical models, those need to be verified and validated against available earthquake recordings, cf. e.g. \citep{bielak2010shakeout,  Paolucci2015, Burks2014}. Nowadays, thanks to the availability of openly accessible data, as for example \cite{SRCMOD, USGS}, this challenge can be tackled in specific regions of the world. PBS generates synthetic time histories of displacement, velocity, acceleration and also other engineering relevant quantities, such as strains, stresses and rotations. The numerical methods used are most often finite differences \cite{Chaljub2015, Pitarka2021, McCallen2021}, finite elements \cite{Bielak2005}, finite volumes \cite{dumbser2007arbitrary, pelties2012three, duru2020stable, breuer2017edge} or, as used here in combination with a linear visco-elastic model, spectral elements in conforming \cite{komatitsch2004simulations, stupazzini2009near, Galvez2014} but also discontinuous ways \cite{DeBasabe2008, Antonietti2018, ferroni2016dispersion}.\\

In this work, we simulate within a single computational model a full seismic event, from source-to-site, and study the effects of ground shaking on a larger building structure. We therefore employ a, mathematically general seismic wave-popagation model to a computational domain consisting of several layers of soil, each with its own material properties and with lengths up to the 100\,km scale. On the small scale of 10-100\,m, with the same mathematical model, we consider a, comparably small, dam structure in order to analyze its performance under the seismic impact of an earthquake. Together this results in a large, multi-scale problem. The challenge of coupling the multiple, non-conforming meshes of different sizes \cite{melas2021three} is tackled by means of the discontinuous Galerkin spectral element method \cite{antonietti2012non}. In addition, two more important factors are considered for properly simulating a dam subjected to seismic excitation: (i) the structure is embedded into its surrounding and therefore the topography should be accurately described, resulting in complex geometries and hence mesh structures, (ii) in contrast to free-standing typical edifices \cite{Building_City} where the surrounding air is most often ignored, the seismic behaviour of the reservoir lake located behind the dam cannot be neglected, and therefore a coupled elasto-acoustic wave-propagation problem needs to be solved. These two factors are considered in this work by making use of digital elevation maps, that can be obtained freely from \cite{GEBCO, SRTM}, to obtain a realistic topographic profile and second by resolving the reservoir lake behind the dam as well. As a matter of fact, within the reservoir lake the propagation of acoustic waves will be modeled by a scalar wave equation. The elasto-acoustic dam-water- and ground-water-interfaces are equipped with force exchange coupling conditions. For a mathematical discussion of the coupled problem, we refer to \cite{MOXElastoAcousticCoupling, flemisch2006elasto}, while \cite{muhr2021hybrid} considers the problem even in a nonlinear, acoustic context. As a seismic source we consider a kinematic rupture model with a prescribed slip-vector and moment-tensor distribution across a fault plane \cite{FaccioliQuarteroni}. The model is also compared to another sourcing mechanism using a plane wave input of a recorded seismogram.\\

The mathematical model and the numerical simulation are validated with respect to the specific magnitude 7 seismic event that took place on the 30$^{\rm th}$ of October in 2020 at around 11:51 h. Its hypocenter lies at 37.8973$^{\circ}\,$N, 26.7953$^{\circ}\,$E in the Icarian Sea northern the isle of Samos, Greece. Approximately 30\,km north-east of it on Turkish mainland there lies the Tahtal{\i}-dam with its fresh-water reservoir. Due to its proximity to the source, the dam was severely threatened by the seismic event; however, the reconnaissance team provided the evidence that no severe damages occurred \cite{Report_EERI}. Besides the topography data mentioned above, in this work, we make use of the large amount of available data regarding the seismic fault source \cite{USGS} for a realistic simulation of the earthquake's origin, recorded seismograms made freely available by the Turkish Disaster \& Emergency Management Authority AFAD \cite{AFAD} and ground material data \cite{USGS}, in order to validate our model and to yield a realistic description of and reliable results for the considered case study.\\

Our simulations have been obtained using the code \textit{SPEED} \cite{mazzieri2013speed, antonietti2012non, stupazzini2009near} employing hexahedral meshes with higher order spectral elements. Sub-meshes with non-matching grids, e.g., at material interfaces with different refinements are coupled using a discontinuous Galerkin approach. Real case simulations as in the present case can easily result in millions of degrees of freedom in space and time, especially when considering higher polynomial orders. Therefore \textit{SPEED} employs a hybrid MPI/OpenMP parallel implementation allowing to harness a large amount of computational resources.\\

We organize the rest of the paper as follows. Section \ref{sec:MathematicalModel} introduces the elasto-acoustic mathematical model, the equations, boundary- and coupling conditions used. In Section \ref{sec:NumericalMethods}, we discuss the adopted discretization in space and time. We derive the semi-discrete form of the model equations by means of spectral elements and introduce a standard time integration scheme.
In Section \ref{sec:SeismicScenario}, we briefly describe the seismic event adopted as a case study in this work. Section \ref{sec:Meshing}, is then devoted to the geometry acquisition (topography/mechanical properties) from real data, and some comments on the mesh generation are given. Finally, the numerical simulations and results are discussed in Section \ref{sec:NumericalSimulations}, where different simulation methods are presented, validated and finally compared.



\section{Mathematical model}\label{sec:MathematicalModel} We begin by defining the mathematical models used to describe the seismic problem, being the elastic model for the solid parts (soil layers and dam in the specific case study) and acoustic model equations for the fluid part (reservoir lake) each with their corresponding sets of boundary, initial and interface conditions. The computational domain together with its individual material subdomains is then also introduced where, without loss of generality, we refer to the specific case study of the Tahtal{\i}-dam considered in this work.  For the whole manuscript, we will denote scalar quantities by regular, greek or latin characters, vectorial quantities will be bold and underlined and tensorial quantities will just be bold.\\

\paragraph{\bf Mathematical model} As a mathematical model for the description of the individual solid parts/subdomains $\Omega_{\textup{e},i},i=1,2,\dots,N_{\textup{e}}$ of the problem, we use the equations of displacement-based linear elasticity \eqref{eq:LinearElasticity} subdomain-wise with Hook's law $\tens{\sigma} = \lambda\textup{tr}(\tens{\epsilon})\mathds{1}+2\mu\tens{\epsilon}$ as constitutive relation \cite{lorenzo2018numerical,Building_City}, $\tens{\epsilon}=\frac{1}{2}\left(\nabla\vecc{u}+\nabla\vecc{u}^{\top}\right)$ being the symmetric gradient of the displacement $\vecc{u}$, and $\lambda$ and $\mu$ being subdomain-wise constant Lamé-parameters, reading
\begin{align*}
	\lambda= \lambda(\vecc{x})=	\lambda_i, ~~\textup{for }\vecc{x}\in\Omega_{\textup{e},i} \textup{ and } i=1,2,\dots,N_{\textup{e}},\\
	\mu= \mu(\vecc{x})=\mu_i, ~~\textup{for }\vecc{x}\in\Omega_{\textup{e},i} \textup{ and } i=1,2,\dots,N_{\textup{e}}.
\end{align*}

\begin{figure}[h!]
	\begin{center}
		\includegraphics[trim=10cm 2cm 9cm 5cm, clip, scale=0.175]{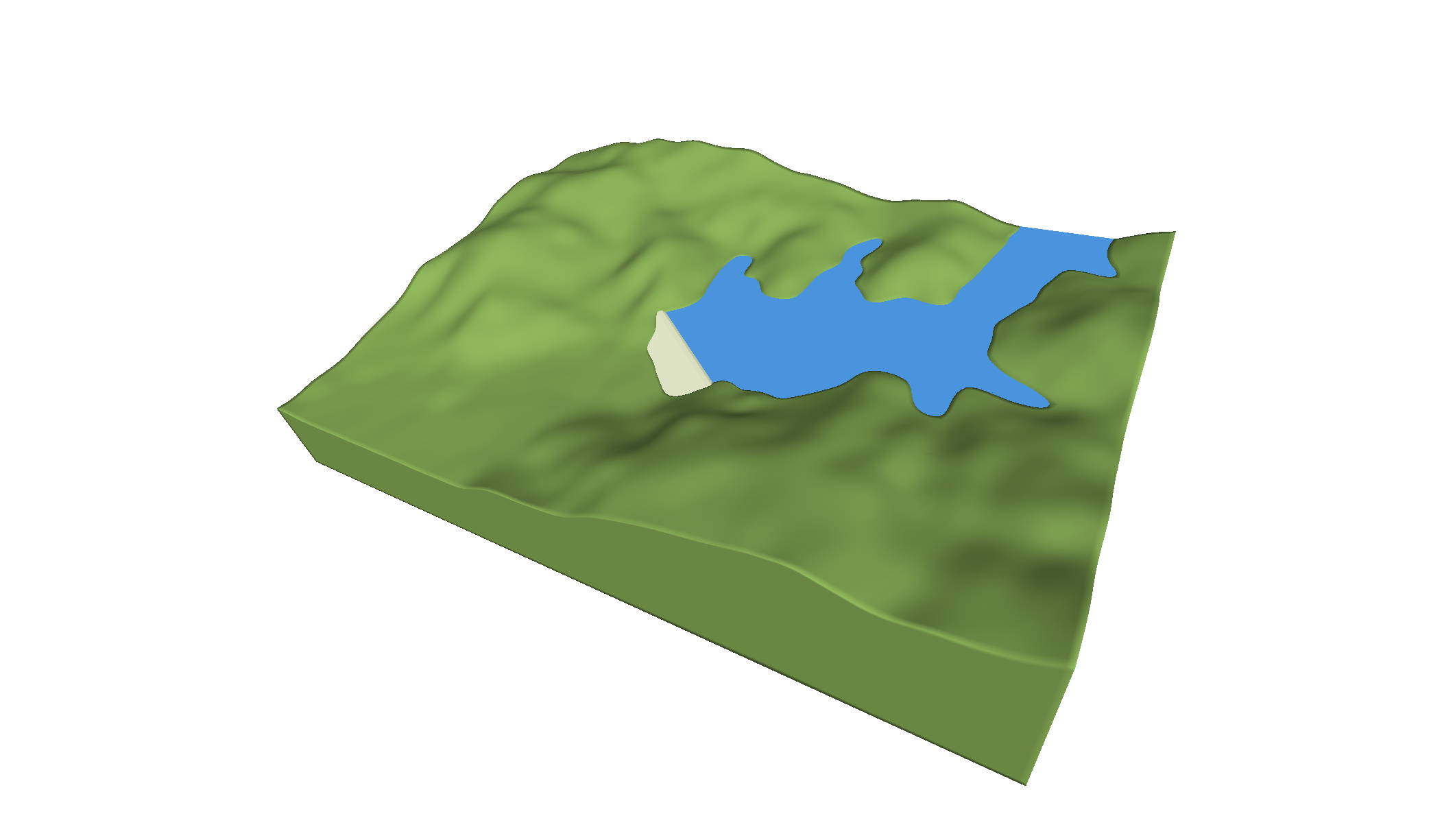}\includegraphics[trim=10cm 2cm 9cm 5cm, clip, scale=0.175]{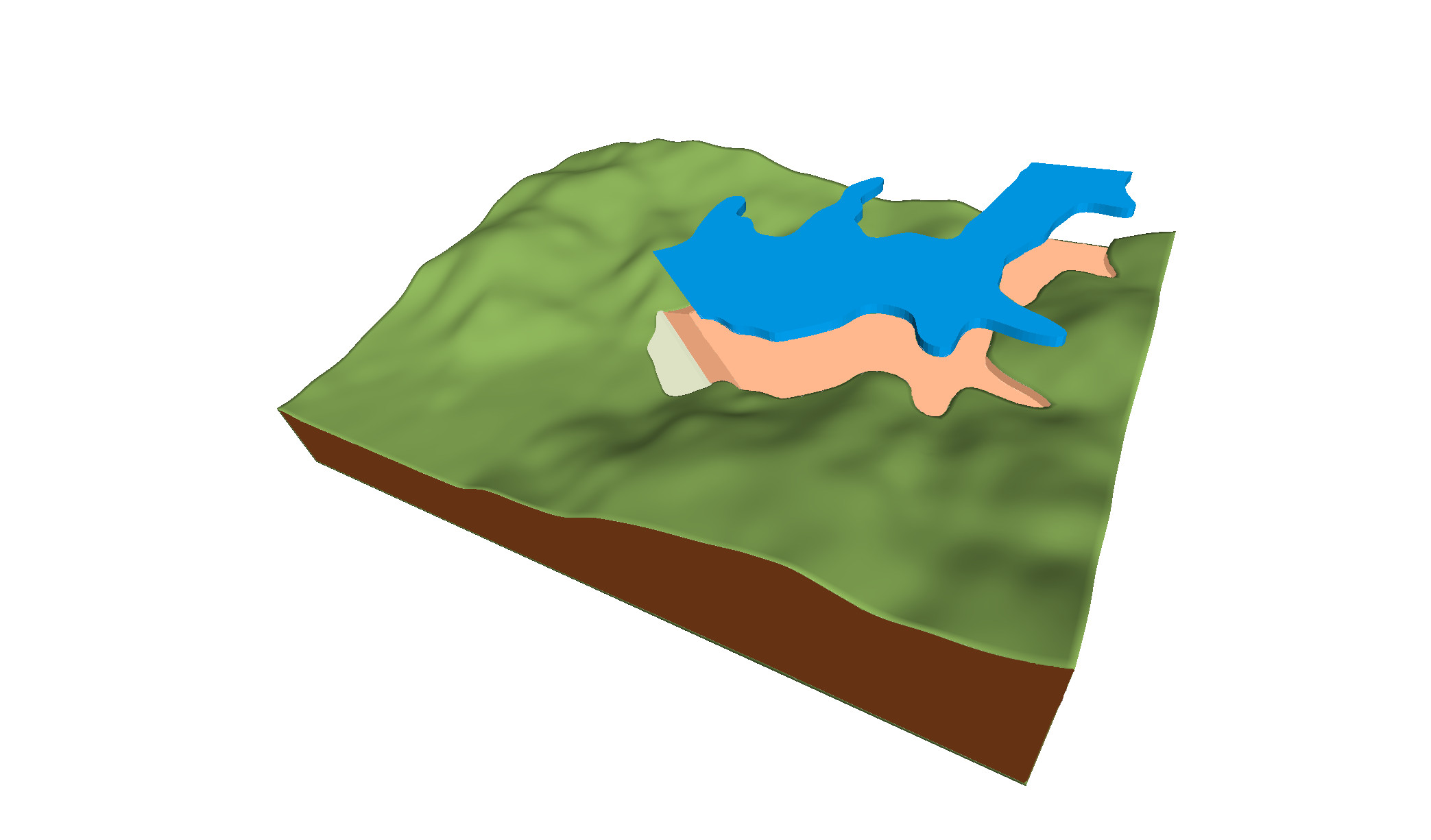}
		\caption{\footnotesize \textbf{(left)} \textrm{Computational domain with three subdomains.} \textbf{Green:} Solid ground $\Omega_{\textup{e},1}$, \textbf{White:} Solid dam structure $\Omega_{\textup{e},2}$, \textbf{Blue:} Acoustic water domain $\Omega_{\textup{a}}$. \textbf{(right)} \textrm{Highlighting boundary/interface conditions}: \textbf{Brown:} Absorbing boundary conditions (all 4 sides, incl. the water cut-off surface plus bottom surface), \textbf{Orange:} Elasto-Acoustic coupling interface.\label{Fig:Domain}, \textbf{Green, White, Water top surface}: Free surfaces, \textbf{Non visible:} Interfaces between $\Omega_{\textup{e},1}$ and $\Omega_{\textup{e},2}$; they are analogous to the elasto-acoustic interfaces.}
	\end{center}
\end{figure}

Parts of $\partial\Omega_{\textup{e}}$ on the top surface with \emph{no} overlying body of water (green and white visible surfaces in Fig.\,\ref{Fig:Domain}, right), summarized as $\Gamma_{\textup{e,N}}$ are treated as free surfaces \eqref{eq:FreeSurface}, the four artificial boundaries in $x$ and $y$ directions as well as the plane bottom surface in $z$-direction (brown in Fig.\,\ref{Fig:Domain}) are equipped with absorbing boundary conditions \eqref{eq:ElasticABCs}. Herein $\vecc{t}^{\ast}$ is a fictitious traction force reducing the amount of artificial reflections originating from these surfaces \cite{stacey1988improved, FaccioliQuarteroni,antonietti2012non}. Parts of $\partial\Omega_{\textup{e}}$ that are interfaces to the acoustic domain $\Omega_{\textup{a}}$, denoted by $\Gamma_{\textup{EA}}$ (orange in Fig.\,\ref{Fig:Domain}), are equipped with non-homogeneous Neumann conditions \eqref{eq:ElasticEA} acting as force-exchange interface conditions to the acoustic field \cite{MOXElastoAcousticCoupling, muhr2021hybrid}. Here the short hand notation of
\begin{equation*}
	\tilde{\psi}:=\psi+\frac{b}{c^2}\dot{\psi}
\end{equation*}
is introduced. On internal interfaces between the individual elastic sub-domains $\Omega_{\textup{e},i}$, collectively denoted by $\Gamma_{\textup{DG}}$ (not visible in Fig.\,\ref{Fig:Domain}, however analogously to elasto-acoustic-interfaces but between ground and dam subdomain), transmission conditions \eqref{eq:ElasticDGCond} are employed, where
\begin{equation*}
	\jump{\tens{\sigma}} := (\tens{\sigma}^+-\tens{\sigma}^-)\vecc{n},\qquad \qquad \jump{\vecc{u}} := (\vecc{u}^+-\vecc{u}^-)\otimes\vecc{n}.
\end{equation*}
Hereby $\vecc{n}$ is the interface normal with arbitrary but fixed orientation and $\iota^{\pm}(\vecc{x}):=\lim_{t\downarrow 0}\iota(\vecc{x}\pm t\vecc{n}),\iota\in\lbrace\tens{\sigma},\vecc{u}\rbrace$. Finally suitable initial conditions \eqref{eq:ElasticIC} for displacement and velocity $\vecc{v}=\dot{\vecc{u}}$ are prescribed, completing the elastic problem.

\begin{alignat}{2}\label{eq:LinearElasticity}
	\rho_{\textup{e},i}\left(\ddot{\vecc{u}}+2\zeta_i\dot{\vecc{u}}+\zeta_i^2\vecc{u}\right)-\nabla\cdot\tens{\sigma(\vecc{u})}&=\vecc{f}, \qquad &&\textup{in } \Omega_{\textup{e},i}\times (0,T], \\
	\tens{\sigma}\vecc{n}&=\vecc{0}, \qquad &&\textup{on } \Gamma_{\textup{e,N}}\times (0,T], \label{eq:FreeSurface}\\
	\tens{\sigma}\vecc{n}&=\vecc{t}^{\ast}, \qquad &&\textup{on } \Gamma_{\textup{e,ABC}}\times (0,T], \label{eq:ElasticABCs}\\
	\tens{\sigma}\vecc{n} &= -\rho_{\textup{a}}\dot{\tilde{\psi}}\vecc{n}, \qquad &&\textup{on } \Gamma_{\textup{EA}}\times (0,T], \label{eq:ElasticEA}\\
	\jump{\tens{\sigma}}=\vecc{0}, \jump{\vecc{u}}&=\tens{0}, \qquad &&\textup{on } \Gamma_{\textup{DG}}\times (0,T], \label{eq:ElasticDGCond}\\
	(\vecc{u},\dot{\vecc{u}}) &= (\vecc{u}_0,\vecc{u}_1), \qquad &&\textup{at }\Omega_{\textup{e}}\times\lbrace 0 \rbrace. \label{eq:ElasticIC}
\end{alignat}

In the above system of equations $\rho_{\textup{e},i}$ are the mass densities of the subdomains $\Omega_{\textup{e},i},i=1,2,\dots,N_{\textup{e}}$,  $\rho_{\textup{a}}$ is the mass density for the acoustics domain and $\zeta_i$, $i=1,\dots,N_{\textup{e}}$ are viscous damping factors proportional to the inverse of time. For future use, we also introduce the compressional $v_p$ and shear $v_s $wave velocities defined as $v_{p,i} = \sqrt{(\lambda_i + 2 \mu_i)/\rho_{\textup{e},i}}$ and $v_{s,i} = \sqrt{\mu_i/\rho_{\textup{e},i}}$ for $i=1,\dots,N_{\textup{e}}$, respectively.\\

In the acoustic subdomain $\Omega_{\textup{a}}$, the linear, damped wave equation \eqref{eq:WaveEquation} in potential form, $\psi$ being the acoustic potential, is used to model the propagation of pressure waves with speed of sound $c$ and damping coefficient $b$. On free water-surfaces (Fig.\,\ref{Fig:Domain} left, blue) homogeneous Neumann conditions \eqref{eq:FreeSurfaceAc}, on artificially generated surfaces, resulting from the cut-out of $\Omega$ from the Earth (Fig.\,\ref{Fig:Domain} backside, where the lake is cut-off), absorbing boundary conditions \eqref{eq:AcousticABCs} \cite{engquist1977absorbing,shevchenko2012self} and on interfaces with the elastic bodies $\Omega_{\textup{e},i}$ (Fig.\,\ref{Fig:Domain} right, orange) interface conditions \eqref{eq:AcousticEA} once more as in \cite{MOXElastoAcousticCoupling, muhr2021hybrid} are imposed. Again, suitable initial conditions \eqref{eq:AcousticIC} complete the acoustic problem.
\begin{alignat}{2}\label{eq:WaveEquation}
	\frac{1}{c^2}\ddot{\psi}-\Delta\tilde{\psi}&=0, \qquad &&\textup{in } \Omega_{\textup{a},i}\times (0,T], \\
	\nabla\tilde{\psi}\cdot\vecc{n}&=0, \qquad &&\textup{on } \Gamma_{\textup{a,N}}\times (0,T], \label{eq:FreeSurfaceAc}\\
	\nabla\tilde{\psi}\cdot\vecc{n}&=-\frac{1}{c}\dot{\psi}, \qquad &&\textup{on } \Gamma_{\textup{a,ABC}}\times (0,T], \label{eq:AcousticABCs}\\
	\nabla\tilde{\psi}\cdot\vecc{n}&= -\dot{\vecc{u}}\cdot\vecc{n}, \qquad &&\textup{on } \Gamma_{\textup{EA}}\times (0,T], \label{eq:AcousticEA}\\
	(\psi,\dot{\psi}) &= (\psi_0,\psi_1), \qquad &&\textup{at }\Omega_{\textup{a}}\times\lbrace 0 \rbrace \label{eq:AcousticIC}	
\end{alignat}

Note that for all material parameters $\iota\in\lbrace\rho_{\textup{e},i},\rho_{\textup{a}},\zeta_i, \lambda_i,\mu_i,c,b\rbrace$ defined above, we assume the existence of uniformly positive and finite bounds above and below. Note that quantities of interest like acoustic or seismic/elastic pressure can be computed from the solutions of the above models via $p_{\textup{ac}}=\rho_{\textup{a}}\dot{\psi}$, $p_{\textup{el}}=-\frac{1}{3}\sum_{i=1}^3\tens{\sigma}_{ii}$ in a post processing step.

\begin{remark}
	We remark that $\Gamma_{\textup{DG}}$ does not have to contain all internal elastic interfaces necessarily and hence the (discontinuous Galerkin) transmission conditions do not have to be applied to all of them. As an alternative, if grids are matching, also a conforming coupling would be possible.
\end{remark}


\paragraph{\bf Computational domain} For the computational core-domain $\Omega$ in the specific case study we choose an, in $x-y$ direction rectangular, cut-out of the Earth around the location of the dam. In $z$-direction $\Omega$ is limited by a plane surface at a given depth below the Earth surface at bottom, while on top the topographic profile of the Earth is used. In a second step, $\Omega$ is divided into an acoustic part $\Omega_{\textup{a}}$ consisting of a portion of the dam reservoir lake lying within $\Omega$, and the remaining part $\Omega_{\textup{e}}$ consisting of solid ground (mountain range, soil layers) and the dam itself. The elastic subdomain $\Omega_{\textup{e}}$ is then further divided into individual subdomains $\Omega_{\textup{e},i}, i=1,2,\dots,N_{\textup{e}}$, each representing an individual material block, with its own set of constant material parameters, which are possibly discontinuous across the interfaces between the blocks. As a prime example, the dam structure would be one such subdomain while the surrounding ground would be an other. See Fig.\,\ref{Fig:Domain}, left, for a graphical depiction of the domain of interest with $N_{\textup{e}}=2$ elastic subdomains.\\

\section{Numerical methods}\label{sec:NumericalMethods}
This section starts with a description of the spatially discrete setting used to approximate a weak solution to the seismic problem \eqref{eq:LinearElasticity}-\eqref{eq:AcousticIC}. We then derive the semi-discrete equation in variational and matrix-vector form and end with some notes about the used time integration scheme.


\subsection{Spatial discretization} In this subsection, the mesh(es) to be used and the finite element spaces built on them are introduced. We hereby closely follow \cite{muhr2021hybrid}, where more details on mesh assumptions being sufficient to prove convergence in a similar setting are given.\\

\paragraph{\bf Meshing} The meshing of the computational domain $\Omega$ is done sub-domain wise. This means that each $\Omega_{\textup{e},i}$ as well as $\Omega_{\textup{a}}$ is subdivided into a mesh $\mathcal{T}_{\textup{e},i}$, resp. $\mathcal{T}_{\textup{a}}$ of hexahedral elements individually. While we assume that also after meshing the discrete interface manifolds do coincide - seen from both sides of the interfaces - the individual meshes on those manifolds do \emph{not} have to. Hence, on the interface the mesh from one side could be a refinement of the mesh from the other side or, for example at a flat interface, could also be a staggered or even completely different mesh. Such non-conformities will be treated by a discontinuous Galerkin approach for the interface coupling.\\ 

\paragraph{\bf Discrete spaces} Our goal is to approximate a weak solution to the seismic problem by means of spectral finite elements. Within each subdomain $\Omega_{\textup{e},i}$ and $\Omega_{\textup{a}}$, this is done in a conforming way, while only at (part of) the interfaces the aforementioned DG approach will be used. This goes hand in hand with the different material properties of the individual sub-domains, and it also allows to keep the amount of degrees of freedom low within the sub-domains, while being flexible at the interfaces.\\
\indent We denote by $\vecc{V}_h^{\textup{e},i}$ the space of discrete ansatz-functions on the $i$-th elastic subdomain, $i=1,2,\dots,N_{\textup{e}}$, which are elementwise polynomials of order $p_i$ when transformed back to the reference element. By $\vecc{V}_h^{\textup{e}}$ we denote the space of global elastic ansatz-functions, which, restricted to any of the subdomains, are within the subdomain's ansatz-space. Similar $V_h^{\textup{a}}$ denotes the ansatz-space for the acoustic subdomain. The setting hence directly corresponds to the one in \cite{muhr2021hybrid} where the overall coupled elasto-acoustic problem, even though in a medical ultrasound setting, including additional non-linear acoustic terms, was analyzed regarding stability and convergence w.r.t. spatial refinement. Therein also further mathematical details are given in a similar setting.\\

\paragraph{\bf Semi-discrete form} With introducing the mean operator defined as:
\begin{equation*}
\mean{\tens{\sigma}} = \frac{1}{2}(\tens{\sigma}^++\tens{\sigma}^-),\qquad \qquad \mean{\vecc{u}} = \frac{1}{2}(\vecc{u}^++\vecc{u}^-).
\end{equation*}
in addition to the already defined jump operator, the semi-discrete variational problem is given by:\\

For any time  $t\in(0,T]$ find $(\vecc{u}_h,\psi_h)\in\vecc{V}_h^{\textup{e}}\times V_h^{\textup{a}}$ such that for all $(\vecc{w}_h,\phi_h)\in\vecc{V}_h^{\textup{e}}\times V_h^{\textup{a}}$ there holds:
\begin{align}\label{eq:WeakForm}
\sum_{i=1}^{N_{\textup{e}}}\left[(\rho_{\textup{e},i}\ddot{\vecc{u}}_h,\vecc{w}_h)_{\Omega_{\textup{e},i}} + (\rho_{\textup{e},i}2\zeta_i\dot{\vecc{u}}_h,\vecc{w}_h)_{\Omega_{\textup{e},i}} + (\rho_{\textup{e},i}\zeta_i^2\vecc{u}_h,\vecc{w}_h)_{\Omega_{\textup{e},i}} + (\tens{\sigma}(\vecc{u}_h),\tens{\epsilon}(\vecc{w}_h))_{\Omega_{\textup{e},i}}\right]\notag \\
	- (\vecc{t}_h^{\ast},\vecc{w}_h)_{\Gamma_{\textup{e,ABC}}}-\langle \mean{\tens{\sigma}(\vecc{u}_h)},\jump{\vecc{v}_h}\rangle_{\Gamma_{\textup{DG}}}-\langle \jump{\vecc{u}_h},\mean{\tens{\sigma}(\vecc{v}_h)}\rangle_{\Gamma_{\textup{DG}}}+\langle\chi\jump{\vecc{u}_h},\jump{\vecc{v}_h}\rangle_{\Gamma_{\textup{DG}}}\\
	+\left(\rho_{\textup{a}}\dot{\tilde{\psi}}_h\vecc{n},\vecc{w}_h\right)_{\Gamma_{\textup{EA}}} + (\dot{\vecc{u}}_h\cdot\vecc{n},\phi_h)_{\Gamma_{\textup{EA}}}\notag \\
	+\left(c^{-2}\ddot{\psi}_h,\phi_h\right)_{\Omega_{\textup{a}}}+(\nabla\tilde{\psi}_h,\nabla\phi_h)_{\Omega_{\textup{a}}} +\left(c^{-1}\dot{\psi}_h,\phi_h\right)_{\Gamma_{\textup{a,ABC}}}=(\vecc{f}_h,\vecc{w}_h)_{\Omega_{\textup{e}}} \notag
\end{align}
and $\vecc{u}_h(0)=\vecc{u}_0=\vecc{0},~\dot{\vecc{u}}_h(0)=\vecc{u}_1=\vecc{0}, \psi_h(0)=\psi_0=0$ and $\dot{\psi}_h(0)=\psi_1=0$, which corresponds to an initial state at rest for the solid as well as for the acoustic quantities.\\

Finally, for any face in $\Gamma_{\textup{DG}}$, defined as the intersection between opposite elemental faces, we define the penalty parameter $\chi$ as:
\begin{equation*}
	\left.\chi\right|_F := \{\lambda+2\mu\}_A \frac{p_F^2}{h_{F}}
\end{equation*}
being $F$ a DG face shared by the mesh elements $E^+$ and $E^-$, $p_F=\max\{p^+,p^-\}$, $h_F=\min\{h^+,h^-\}$, $\{ q \}_A$ the harmonic average of the quantity $q$, and $\beta$ a positive real number at our disposal, cf. \cite{antonietti2012non, DeBasabe2008, muhr2021hybrid}.\\

\paragraph{\bf Matrix-Vector form} After representing the discrete trial functions $\vecc{u}_h$ and $\psi_h$ as well as their temporal derivatives in the nodal finite element basis, denoting the coefficient vectors by (abuse of notation) also $\vecc{u}_h$ and $\vecc{\psi}_h$, the following system of ODEs in matrix-vector form is directly obtained from \eqref{eq:WeakForm}:
{\small\begin{align}\notag
	\begin{pmatrix}
	\mat{M}_{\textup{e}}^{(2)} & 0 \\
	0 & \mat{M}_{\textup{a}}^{(2)}
	\end{pmatrix}\begin{pmatrix}
	\ddot{\vecc{u}}_h \\ \ddot{\vecc{\psi}}_h
	\end{pmatrix}
	=&
	-\begin{pmatrix}
	\mat{M}_{\textup{e}}^{(1)} & 0 \\
	0 & 0
	\end{pmatrix}\begin{pmatrix}
	\dot{\vecc{u}}_h \\ \dot{\vecc{\psi}}_h
	\end{pmatrix}
	-
	\begin{pmatrix}
	\mat{M}_{\textup{e}}^{(0)}& 0 \\
	0 & 0
	\end{pmatrix}\begin{pmatrix}
	{\vecc{u}}_h \\ \vecc{\psi}_h
	\end{pmatrix}\\ \label{eq:MatrixEQ}
	&+\begin{pmatrix}
	-\mat{K}_{\textup{e}}+\mat{D}+\mat{D}^{\top}-\mat{P} & 0 \\
	0 & -\mat{K}_{\textup{a}}
	\end{pmatrix}\begin{pmatrix}
	{\vecc{u}}_h \\ \vecc{\psi}_h
	\end{pmatrix}
	\\  \notag 
	&-\begin{pmatrix}
	0 & \mat{E} \\
	\mat{A} & 0
	\end{pmatrix}\begin{pmatrix}
	\dot{\vecc{u}}_h \\ \dot{\vecc{\tilde{\psi}}}_h
	\end{pmatrix}
	-
	\begin{pmatrix}
	0 & 0 \\
	0 &  \mat{C}_{\textup{a}}+\mat{C}_{\textup{a},ABC}
	\end{pmatrix}\begin{pmatrix}
	\dot{\vecc{u}}_h \\ \dot{\vecc{\psi}}_h
	\end{pmatrix}
	+
	\begin{pmatrix}
	\vecc{T}_h^{\ast}\\ 0
	\end{pmatrix}
	+
	\begin{pmatrix}
	\vecc{F}_h\\ 0
	\end{pmatrix}\\[5mm] \notag
	\textup{with initial conditions: }
	&
	\begin{pmatrix}
	\vecc{u}\\ \vecc{\psi}
	\end{pmatrix}=\begin{pmatrix}
	\dot{\vecc{u}}\\ \dot{\vecc{\psi}}
	\end{pmatrix}=\begin{pmatrix}
	\vecc{0} \\ \vec{0}
	\end{pmatrix}
\end{align}}
Herein $\mat{M}_{\textup{e}}^{(1)}$ and $\mat{M}_{\textup{e}}^{(0)}$ are subdomain-wise scaled versions of the standard mass-matrix $\mat{M}_{\textup{e}}^{(2)}$ for the elastic part of the problem. Due to the DG approach for the interface coupling, they are block-diagonal with one block per subdomain. In fact, due to the use of spectral elements, they even become diagonal. This property is not only advantageous to be exploited by the time stepping scheme but also allows to denote the aforementioned scaling using the diagonal matrix $\mat{Z}$, which has entries $\zeta_i$ for each degree of freedom belonging to the elastic subdomain $\Omega_{\textup{e},i}$. With that, the scaled mass matrices read $\mat{M}_{\textup{e}}^{(1)} = 2\mat{Z}\mat{M}_{\textup{e}}^{(2)}$ and $\mat{M}_{\textup{e}}^{(0)} = \mat{Z}^2\mat{M}_{\textup{e}}^{(2)}$. Analogously $\mat{M}_{\textup{a}}^{(2)}$ is the standard mass-matrix for the acoustic part. $\mat{K}_{\textup{a}}$ resp. $\mat{K}_{\textup{e}}$ denote stiffness-matrices for the scalar resp. vectorial case, $\mat{C}_a$ incorporates the $b$-damping of the acoustic equation. On the right hand side $\mat{D}$ and $\mat{P}$ are the matrices stemming from the inter-elastic DG coupling, where $\mat{P}$ contains the penalization term, while $\mat{E}$ and $\mat{A}$ contain the coupling terms between the elastic and acoustic domains. Finally $\mat{C}_{\textup{a},ABC}$ is a boundary mass-matrix used in the implementation of the acoustic absorbing boundary conditions, $\vecc{T}_h^{\ast}$ incorporates the elastic absorbing boundary conditions and $\vecc{F}_h$ the elastic source term inducing the earthquake (cf. Sec.\,\ref{sec:NumericalSimulations} for different source models).


\subsection{Time integration}
The basis scheme used for time integration is the Leap-Frog scheme in its full-step predictor-corrector form \cite{hairer2003geometric}. After the prediction step 
\begin{equation*}
	\vecc{\phi}^{(n+1)} = \vecc{\phi}_{\textup{pred}} = \vecc{\phi}^{(n)} + \Delta t \dot{\vecc{\phi}}^{(n)} + \frac{1}{2}\Delta t^2 \ddot{\vecc{\phi}}^{(n)},\qquad
	\dot{\vecc{\phi}}_{\textup{pred}} = \dot{\vecc{\phi}}^{(n)} + \frac{1}{2}\Delta t \ddot{\vecc{\phi}}^{(n)}
\end{equation*}
for $\vecc{\phi}\in\lbrace\vecc{u}_h,\vecc{\psi}_h\rbrace$ the right hand side of \eqref{eq:MatrixEQ} can be evaluated using the predicted values (and the previous timestep value for $\ddot{\vecc{\psi}}_h$ within $\dot{\tilde{\vecc{\psi}}}_h$). One then conducts a solver step for the left hand side variables $\ddot{\vecc{u}}_h^{(n+1)}$ and $\ddot{\vecc{\psi}}_h^{(n+1)}$ followed by the correction step
\begin{equation*}
	\dot{\vecc{\phi}}^{(n+1)} = \dot{\vecc{\phi}}_{\textup{pred}} + \frac{1}{2}\Delta t  \ddot{\vecc{\phi}}^{(n+1)}
\end{equation*}
again $\vecc{\phi}\in\lbrace\vecc{u}_h,\vecc{\psi}_h\rbrace$. \\

To take into account the multiscale nature in time, we also employ a simple local time-stepping approach for the acoustic part of the problem similar to \cite{shevchenko2012multi}. Especially in the vicinity of the inclined dam surface that is under water (cf. Fig.\,\ref{Fig:DamCAD}, right) water elements become quite thin hence requiring small timesteps that are not necessary in the remaining region of the problem, especially the huge land masses $\Omega_{\textup{e},i}$ of seismic propagation. Hence, within each regular time step of size $\Delta t$ the solution of the acoustic equation within $\Omega_{\textup{a}}$ consisting of prediction, evaluation, solving and correction is internally repeated within a loop of $n_{\textup{loc}}$ timesteps of size $\frac{\Delta t}{n_{\textup{loc}}}$ before the next regular step for the elastic portion of the problem is conducted. In the presented simulation cases a value of $n_{\textup{loc}}=10$ yields good results.


\section{Seismic scenario description}\label{sec:SeismicScenario}
We now provide some selected information about the specific, real seismic event that will serve as a case study for the mathematical model and numerical methods described in Sec.\,\ref{sec:MathematicalModel} and \ref{sec:NumericalMethods}. The event was chosen due to the large amount of seismograms recorded and therefore allowing a proper validation of the simulation.\\

The Samos Island (Aegean Sea) earthquake struck at 14:51 local time in Turkey, on 30 October 2020. Severe damages have been observed in some densely populated districts of $\dot{\textup{I}}$zmir (Bayraklı, Bornova, Karşıyaka and Konak), and 118 fatalities have been reported. The present event occurred in the cross‐border region between the eastern Aegean Sea islands and Western Turkey, which is among the most seismically active areas in Eastern Mediterranean and has been the site of devastating earthquakes in both recent and historical times, see \cite{Guidoboni1994}. The fault that ruptured during the mainshock is located offshore the northern coast of Samos Island, and it was previously identified as Kaystrios Fault (see GreDaSS database \cite{Caputo2013TheGD} and the GEM‐Faults database \cite{GEMFaults}). The geometric characteristics of the Kaystrios Fault are: strike in the range 260° to 290°, dip 45° to 70° and rake in the range ‐100° to ‐80°, while the maximum depth of the fault is estimated as 14.5\,km.\\
\indent The mainshock was well recorded by the broadband seismic networks of Greece and Turkey, and 35 records are located within $100$\,km distance from the epicenter \cite{Cauzzi2021, AFAD-TADAS, Report_EERI,Report_METU,yakut2021performance}.
Between October 31st and November 6th a group of engineers assembled by the $\dot{\textup{I}}$zmir Regional Directorate of State Hydraulic Works (DSI) and a reconnaissance team from METU visited dam sites to document the performance of earth-fill and rock-fill dams shaken by the event.
After a detailed inspection of six small to medium size earth-fill and rock-fill dams including the Tahtal{\i} dam, no apparent damage was reported by reconnaissance teams.\\
\indent As reported by \cite{TosunTosun2018}, the Tahtal{\i} dam is a rockfill dam on the Tahtal{\i} River near Gumuldur County in the $\dot{\textup{I}}$zmir Metropolitan Area. It has a $54.4$\,m height from river bed. When the reservoir is at maximum capacity, the facility impounds $306.6$\,hm$^3$ of water in its reservoir. Its construction was finished in 1999. It was designed to provide domestic water with an active volume of $287$\,hm$^3$.
As reported by \cite{TosunTosun2018}, the Tahtal{\i} dam is only $1.9\,$km away from an active fault and, according to the seismic hazard analyses performed, it will be subjected to a peak ground acceleration of $0.277~g$ by an earthquake of $M_w~5.7$. The Tahtal{\i} dam is at second place, after Gordes dam, when regarding the total capacity of the reservoir and one of the most critical dams in the $\dot{\textup{I}}$zmir Metropolitan Area. 


\section{Geometry acquisition and mesh generation}\label{sec:Meshing}

The methods described in this section are general and can easily be adopted to different scenarios. However, the individual steps are motivated and illustrated based on the case study discussed in Sec.\,\ref{sec:SeismicScenario}.\\

\begin{figure}[h!]
	\begin{center}
		\hspace*{-5mm}\includegraphics[trim=0cm 0cm 0cm 0cm, clip, scale=0.45]{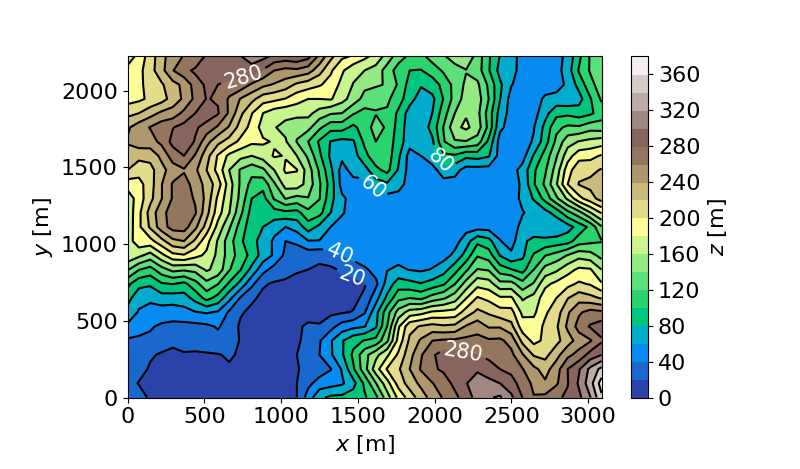}\hspace*{0cm}\raisebox{0.73cm}{\includegraphics[trim=0cm 0cm 0cm 0cm, clip, scale=0.125]{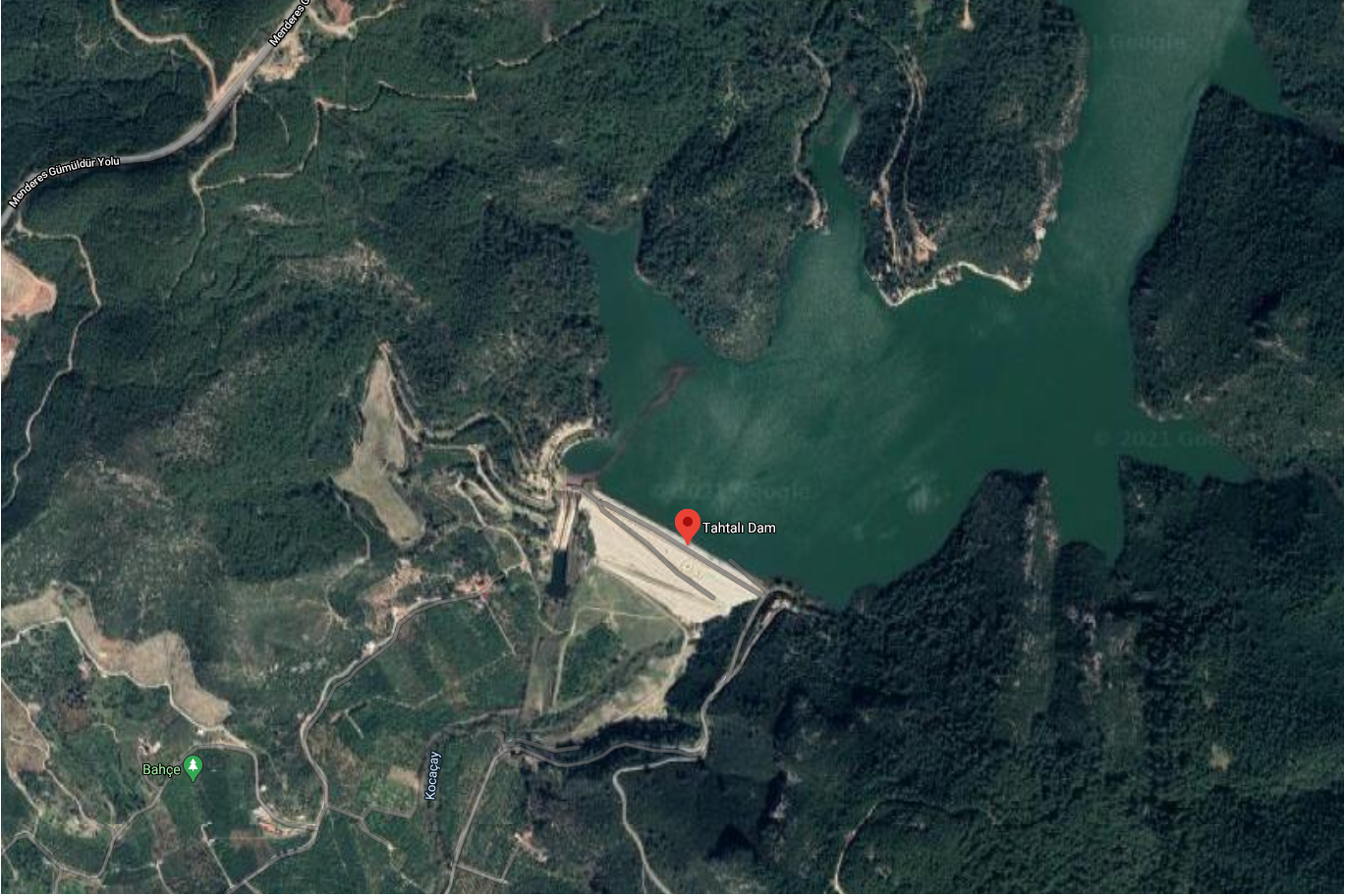}}\\
		\vspace{-7mm}
		\hspace{8.6cm}\parbox{6cm}{\tiny Map data: Google, Imagery \textcopyright 2021 CNES/Airbus, Maxar Technologies, Map data \textcopyright 2021}
		\caption{\footnotesize \textbf{(left)} Topographic map of the region around the Tahtal{\i}-dam with reference coordinate-system. Data source \cite{SRTM}. \textbf{(right)} Satellite view of the same area. The red marker is located at $38.0888\,^{\circ}\textup{N}, 27.0415\,^{\circ}\textup{E}$. Image source \cite{GoogleMaps}.\label{Fig:Topography}}
	\end{center}
\end{figure}

\paragraph{\bf Digital elevation map data} In the presently examined case study, the topography of the surrounding area of the Tahtal{\i}-dam is obtained from the \textit{SRTM Digital Elevation Database} of CGIARCSI \cite{SRTM}, which contains topographic elevation data in a 3 arc second grid.\\
On the left-hand side of Fig.\,\ref{Fig:Topography} a contour plot of a roughly 3\,km$\times$2\,km area around the Tahtal{\i}-dam is shown, while on the right-hand side of Fig.\,\ref{Fig:Topography}, the respective area in a satellite view for a clear image of the dam location is given.\\

\begin{figure}[h!]
	\begin{center}
		\hspace*{0mm}\includegraphics[trim=0cm 0cm 0cm 0cm, clip, scale=0.17]{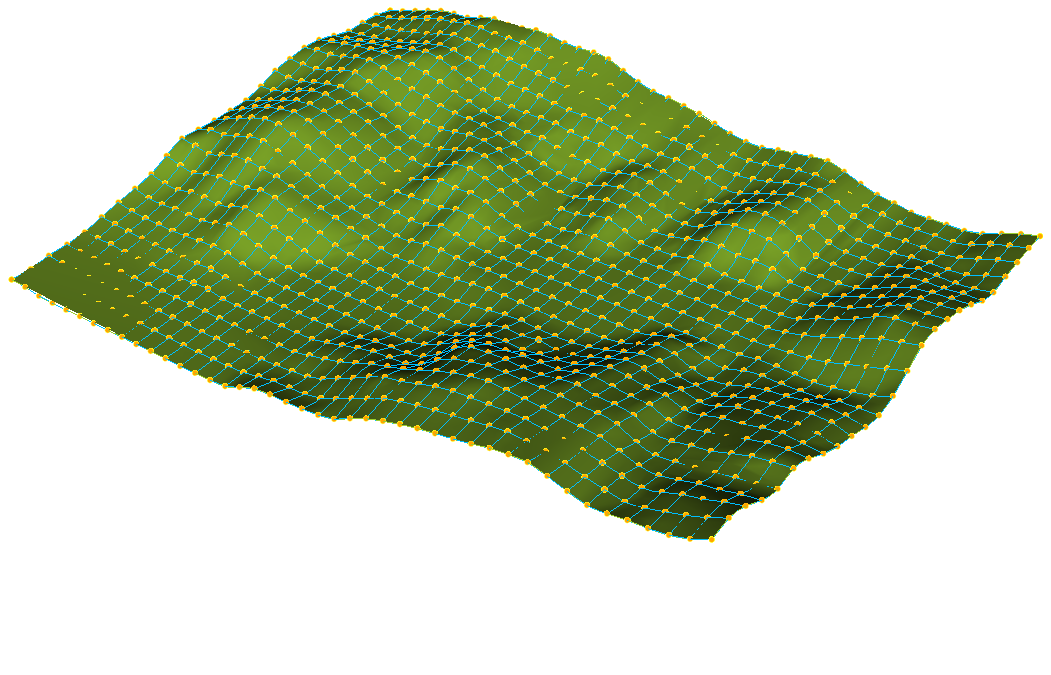}\hspace*{0cm}\raisebox{0cm}{\includegraphics[trim=0cm 0cm 0cm 0cm, clip, scale=0.17]{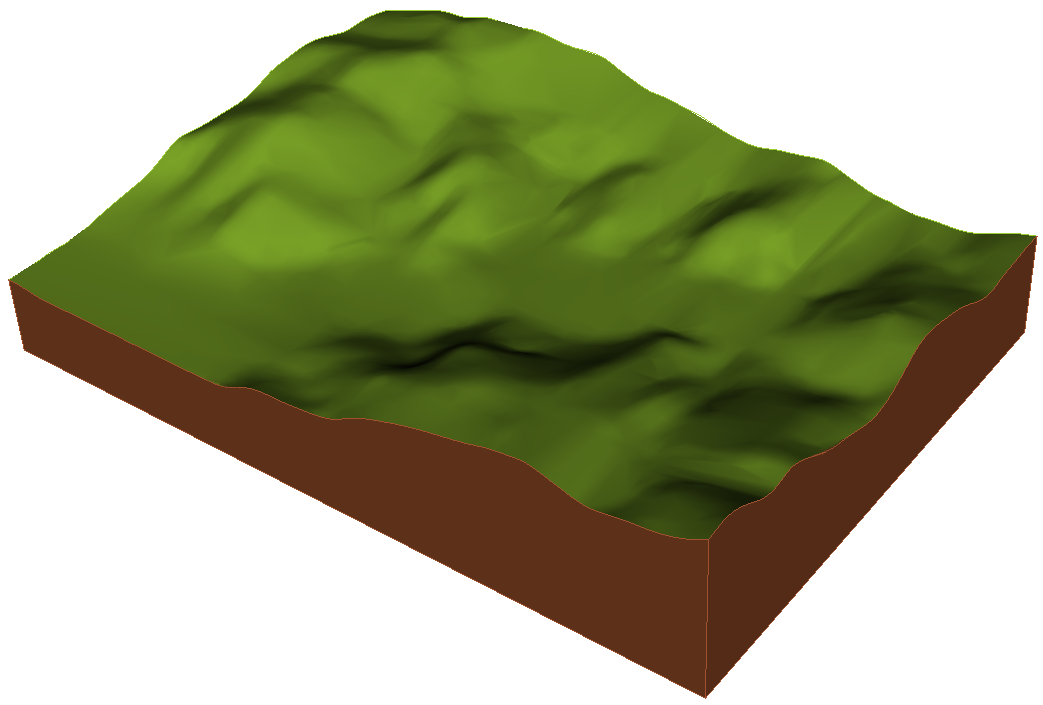}}
		\caption{\footnotesize \textbf{(left)} Orange: Elevation point data obtained from \cite{SRTM}, Blue: $x-y$-spline curves, Green: Spline net-surface approximation of topography. \textbf{(right)} Volumetric ground object with closing sides (and bottom surface at $z=-300\,\textup{m}$) in brown.\label{Fig:GroundReconstruction}}
	\end{center}
\end{figure}

\paragraph{\bf CAD model construction} Starting from the elevation point-cloud data, a volumetric object is constructed and afterwards modified in order to add the dam and the water, and finally proceeding with the meshing. To tackle this challenge the software Cubit \cite{Cubit,casarotti2008cubit}, which also features some CAD-capabilities, was used firstly for the surface- and then volume-reconstruction.\\
Using the point cloud elevation data, a spline-surface approximating the topography is created, see Fig.\,\ref{Fig:GroundReconstruction}. For the bottom surface a rectangle in the $x$-$y$-plane directly below the topographic spline-surface at depth $z=-300\,\textup{m}$ is used. The four sides remaining to close the volume are flat surfaces and connect the four boundary curves of bottom and top surfaces.\\

\begin{figure}[h!]
	\begin{center}
		\hspace*{-15mm}\includegraphics[trim=0cm 0cm 5cm 0cm, clip, scale=0.22]{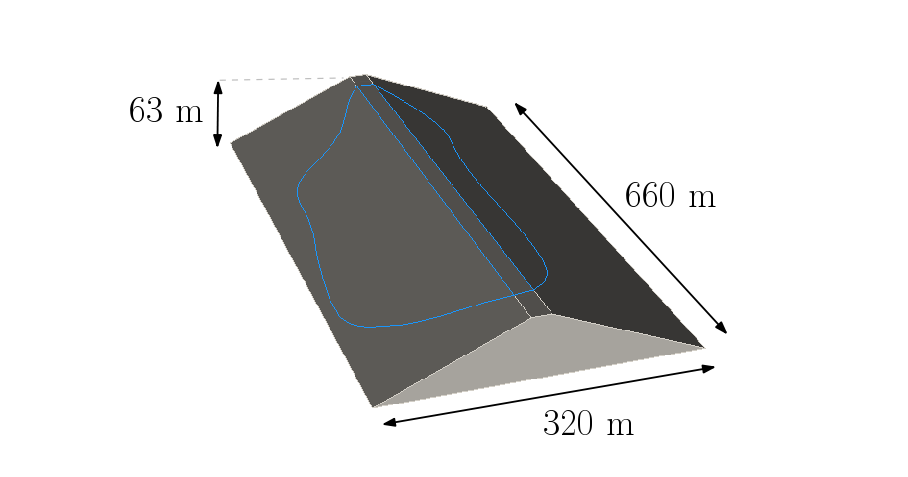}\hspace*{0cm}\raisebox{0cm}{\includegraphics[trim=5cm 0cm 7cm 0cm, clip, scale=0.17]{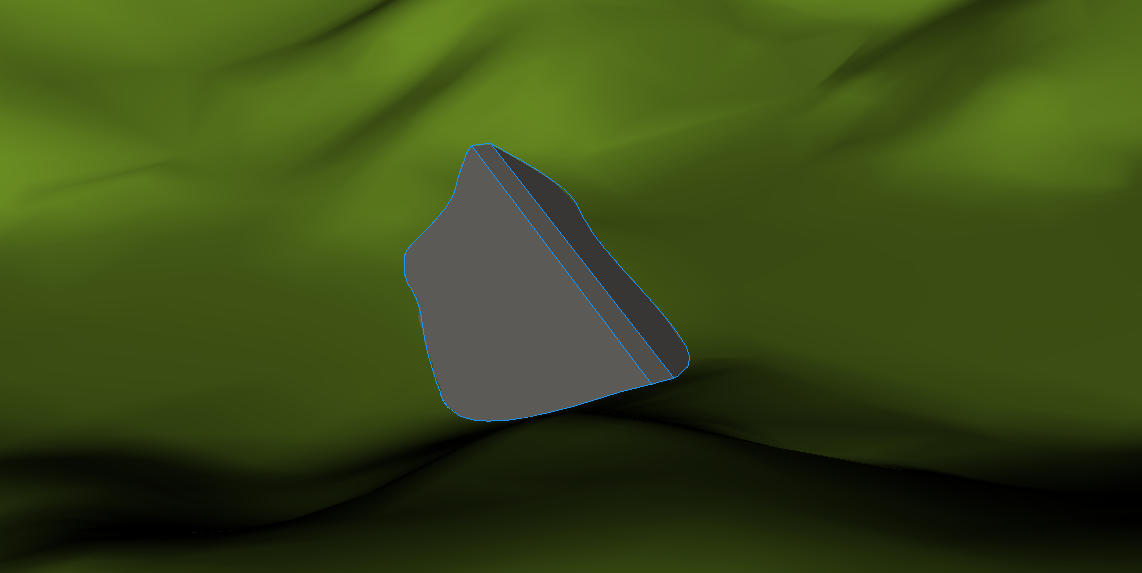}}\hspace*{2mm}\includegraphics[trim=3cm 0cm 5cm 0cm, clip, scale=0.135575]{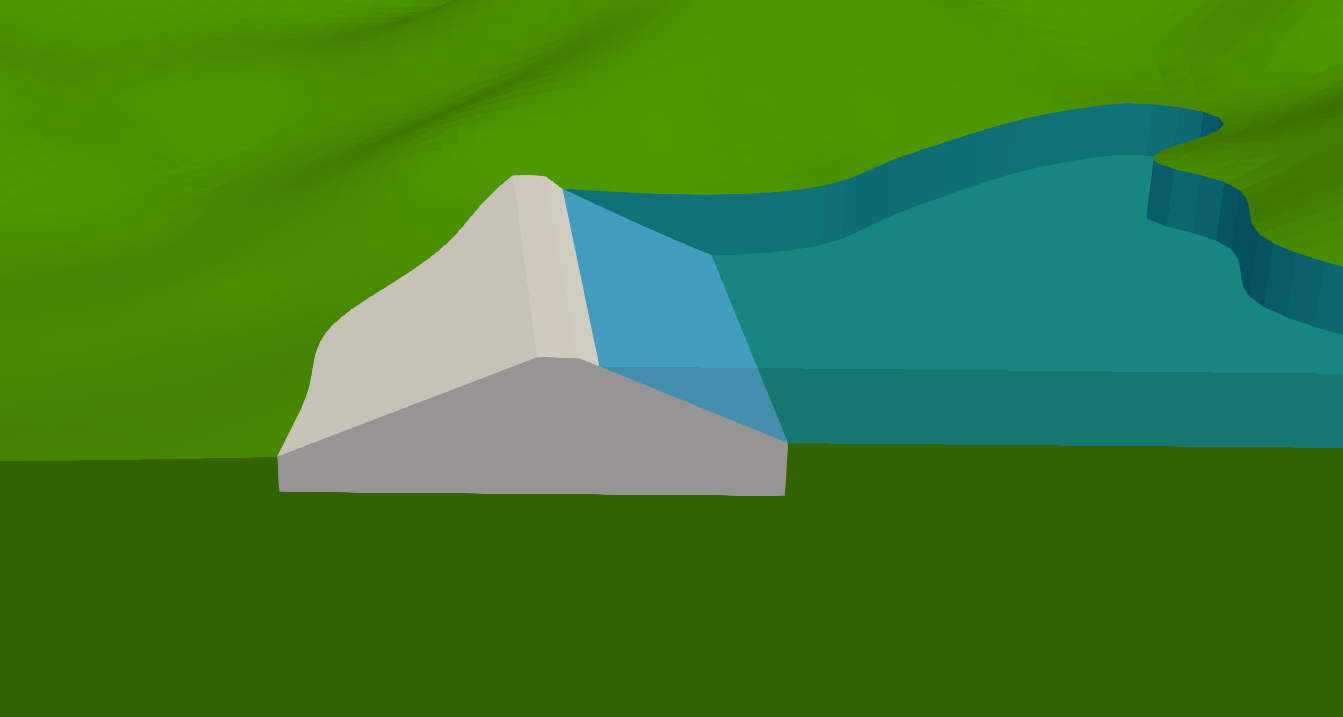}
		\caption{\footnotesize \textbf{(left)} CAD-model of a gravitational dam. The blue lines (compare to figure right) are the intersection curves between the dam and ground model. \textbf{(middle)} Dam-structure embedded into surface topography. \textbf{(right)} Cross section View on CAD-model of ground, dam and water. Due to the embedding into the topography the dam measures reduce to approximately $L\approx 570$\,m, $B\approx 250$\,m and $H\approx 50$\,m of visible size.\label{Fig:DamCAD}}
	\end{center}
\end{figure}

Once the topography is available, the dam is added, see Fig.\,\ref{Fig:DamCAD}, left. The dam is placed in the proper location within the topography and, by means of simple boolean intersection/subtraction as well as downward extrusion operations between the two volumes, is \textit{``naturally''} finally embedded into the surrounding area. Fig.\,\ref{Fig:DamCAD}, middle, shows the portion of the dam after the operations within the topography.\\

Having placed ground and dam, only the acoustic water-subdomain - the reservoir lake - is missing. It is created by ``flooding'' the ground block behind the dam and again some downward extrusion followed by a boolean subtraction from the ground block. The resulting water block can be seen in Fig.\,\ref{Fig:Domain}, right, as a whole, Fig.\,\ref{Fig:DamCAD}, right, gives a cross-sectional view. It should be observed that the inclined dam surface interfacing with the water block is explicitly created as part of the model, again by means of boolean CAD operations. The final model, consisting of the blocks ground, dam and water can be seen in Fig.\,\ref{Fig:Domain}, left.\\

\begin{figure}[h!]
	\begin{center}
		\includegraphics[trim=0cm 0cm 0cm 0cm, clip, scale=0.1008]{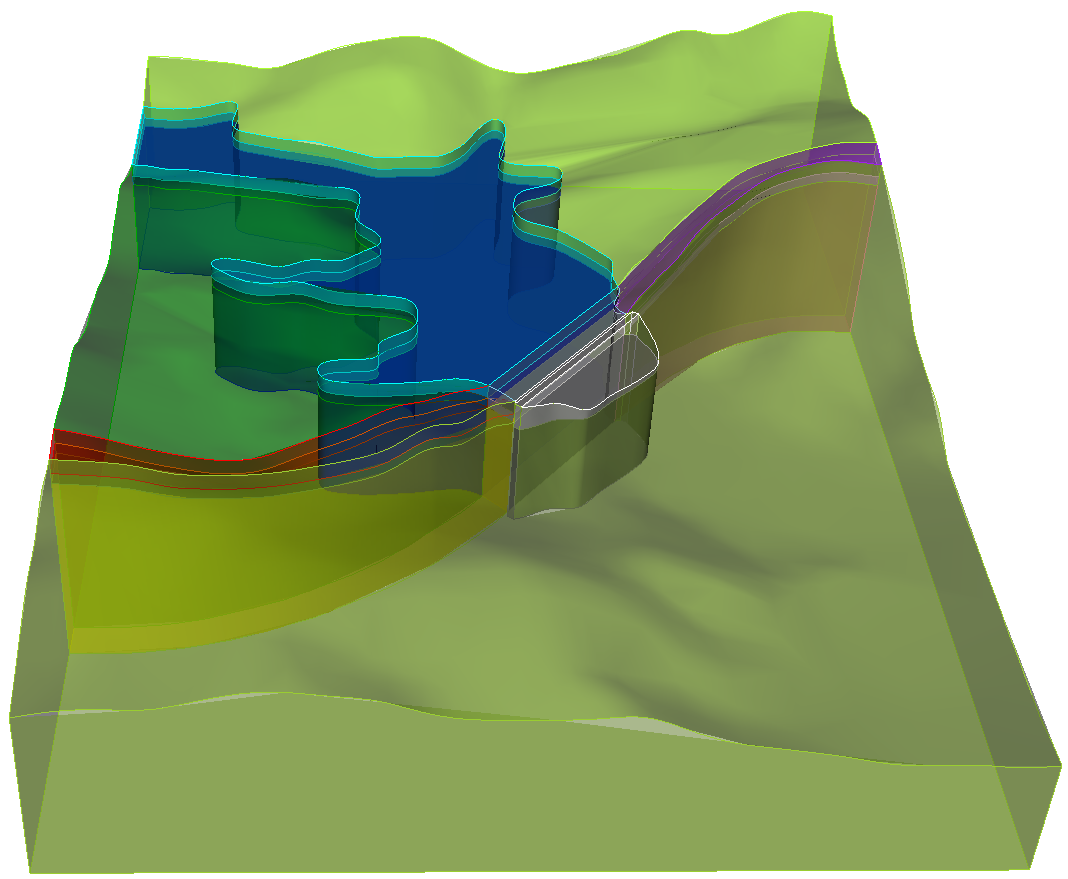}\hspace*{0.5cm}\includegraphics[trim=0cm 0cm 0cm 0cm, clip, scale=0.0864]{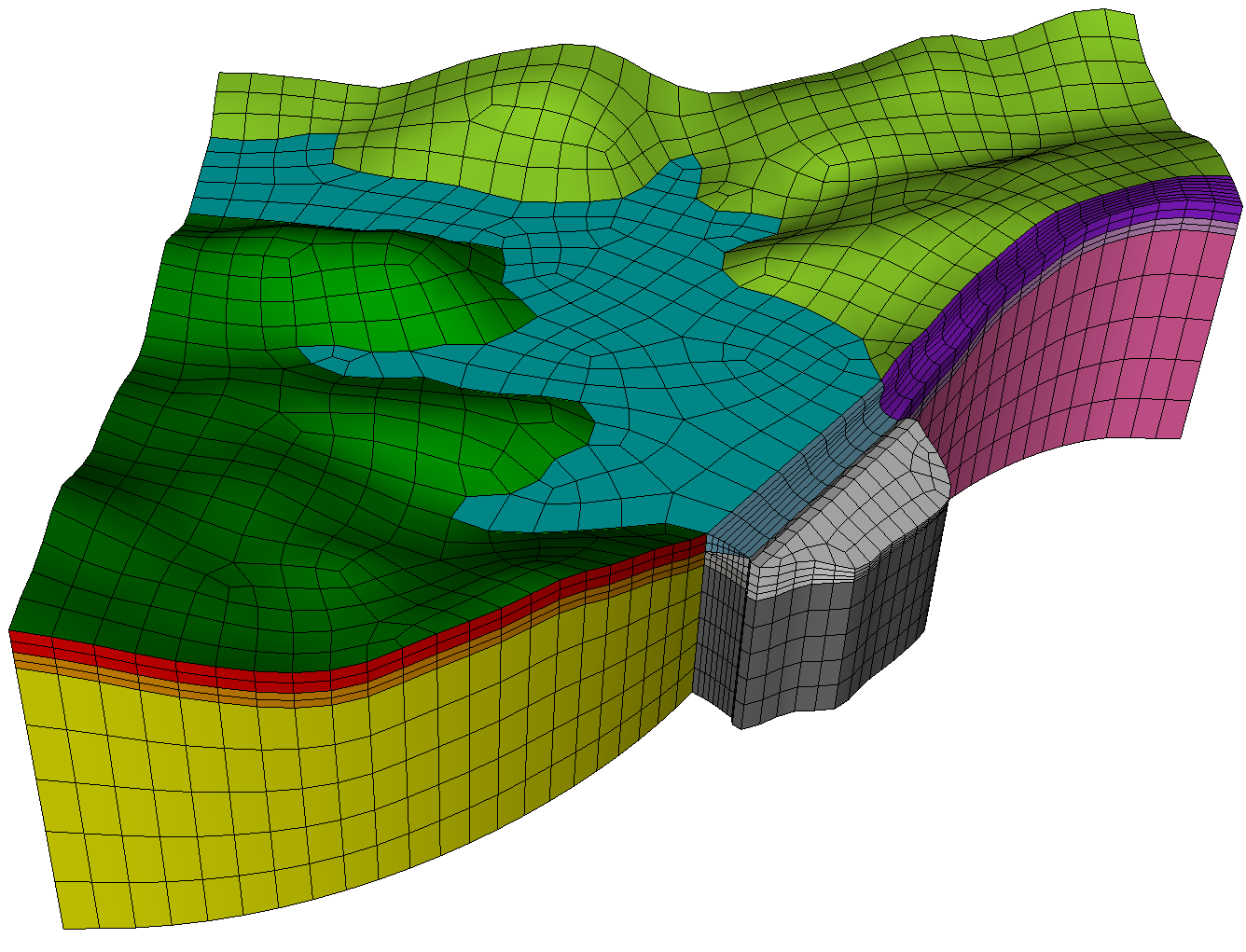}\hspace*{0.5cm}\includegraphics[trim=0cm 0cm 0cm 0cm, clip, scale=0.099]{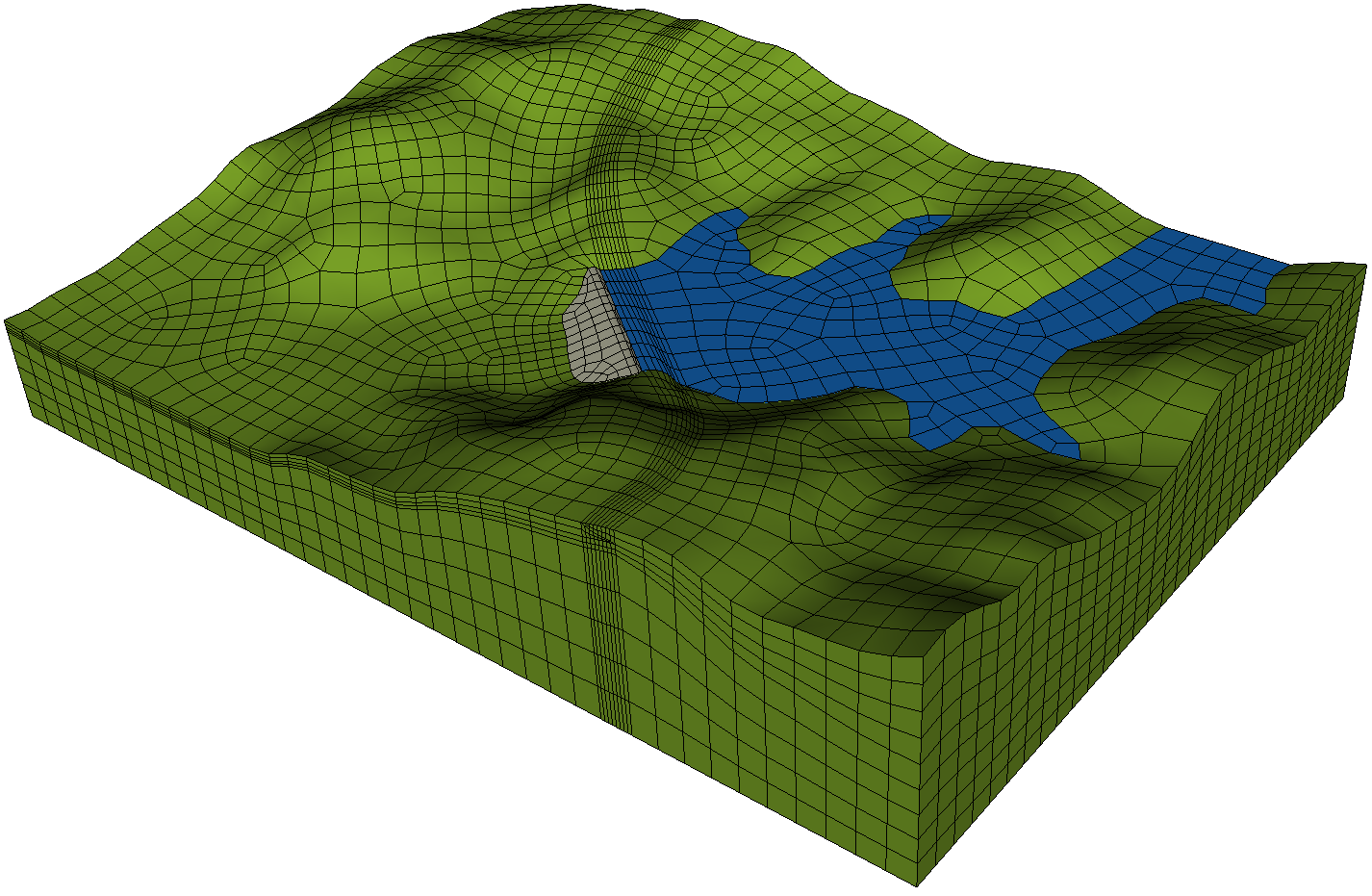}
		\caption{\footnotesize \textbf{(left)} Division into sub-blocks for the \textit{pave-and-sweep} meshing approach. \textbf{(middle)} Meshed sub-blocks with free view on critical water-dam interface. \textbf{(right)} Complete mesh of the computational domain, divided into three blocks: Ground, dam and water. Compare to Fig.\,\ref{Fig:Domain}, left.\label{Fig:Blocking}}
	\end{center}
\end{figure}

\paragraph{\bf Meshing} The usage of a complete hexahedral mesh, with matching interfaces, bears a much more difficult task in mesh generation than a tetrahedral mesh. Even though the software Cubit offers automatic hex-mesh routines like \textit{sculpt} \cite{owen2017hexahedral}, the resulting mesh quality turned out to be not fully satisfactory in the present case, especially in the proximity of interfaces or corners with sharp angles. The process of mesh generation for this work was done in a semi-automatic way, where in a manual pre-processing step the three blocks of the CAD model were subdivided into even more smaller blocks, each with an easily meshable form. Fig.\,\ref{Fig:Blocking}, left, shows all those sub-blocks used.
The meshing of the sub-blocks was then conducted automatically by a \textit{pave-and-sweep} approach \cite{blacker1991paving, blacker1991analysis, mingwu1996multiple}. Special care had to be taken for the inclined water-dam interface (cf. Fig.\,\ref{Fig:DamCAD}, right).\\

After successful meshing of all sub-blocks by the aforementioned strategy, sub-block meshes of the same material were merged again to conforming meshes such that in the end (non-conforming) interfaces only remain, where introduced in Sec.\,\ref{sec:MathematicalModel}. Fig.\,\ref{Fig:Blocking}, right, shows a complete mesh of the three material blocks, cf. Fig.\,\ref{Fig:Domain}, left.


\section{Numerical simulation results}\label{sec:NumericalSimulations}

In this section, we present the numerical results obtained. We conduct essentially three numerical simulations: The first one considers the seismic event on a regional scale and ignores the dam or the reservoir lake but takes into account the topography in order to validate the seismic source and numerical wave-propagation model against observed seismograms. In the second scenario, the vicinity region of the dam/reservoir from Fig.\,\ref{Fig:Domain} and \ref{Fig:Blocking} is excited by a simple plane wave input at the bottom of the domain where the time history of the adopted input signal corresponds to the seismogram recorded in the proximity, namely AFAD \# 3536. Finally, a full source-to-site simulation of the considered event is conducted, spanning the multiple length scales from the seismic fault plane up to the dam structure including the complete domain from the second analysis as a sub-domain.\\
\indent As previously mentioned, the dam and other details are ignored in the regional simulation; this allows to have sufficient degrees of freedom and computational resources available to span a broader region and therefore being capable to incorporate 10 seismograph stations for comparison. In the following source-to-site simulation, the domain is tailored around the optimal size that encompasses the fault plane and the dam region and the thereby saved resources are used to adequately simulate the details in the dam vicinity region.\\
\indent We would like to point out that previously, different verification tests have been considered in order to evaluate the accuracy of the numerical discretization adopted. In particular, we refer the reader to \cite{Antonietti2020_nme} for the linear elasto-acoustic case and to \cite{muhr2021hybrid} for the non-linear and viscous elasto-acoustic case.\\
\indent It is important to recall the concepts of verification and validation \cite{MoczoEtAl2014, MaufroyEtAl2015}: verification of a numerical method may be defined as the demonstration of the consistency of the numerical method with the original mathematical–physical problem defined by the controlling equation, constitutive law, and initial and boundary conditions. The quantitative analysis of accuracy should be a part of the verification. Once the numerical method is analyzed and verified for accuracy, it should be validated using observations. In general, the validation may be defined as the demonstration of the capability of the theoretical model (i.e., the mathematical–physical model and its numerical approximation) to predict and reproduce observations. Normally the criteria and the metrics adopted for the verification and the validation phase are different, given the complexity of the physical problem analyzed.


\subsection{Large scale validation-simulation}\label{subsec:Validation}
\begin{figure}[h!]
	\begin{center}
		\raisebox{3.25cm}{\footnotesize 35\,km}~~~\includegraphics[trim=0cm 0cm 0cm 0cm, clip, scale=0.2]{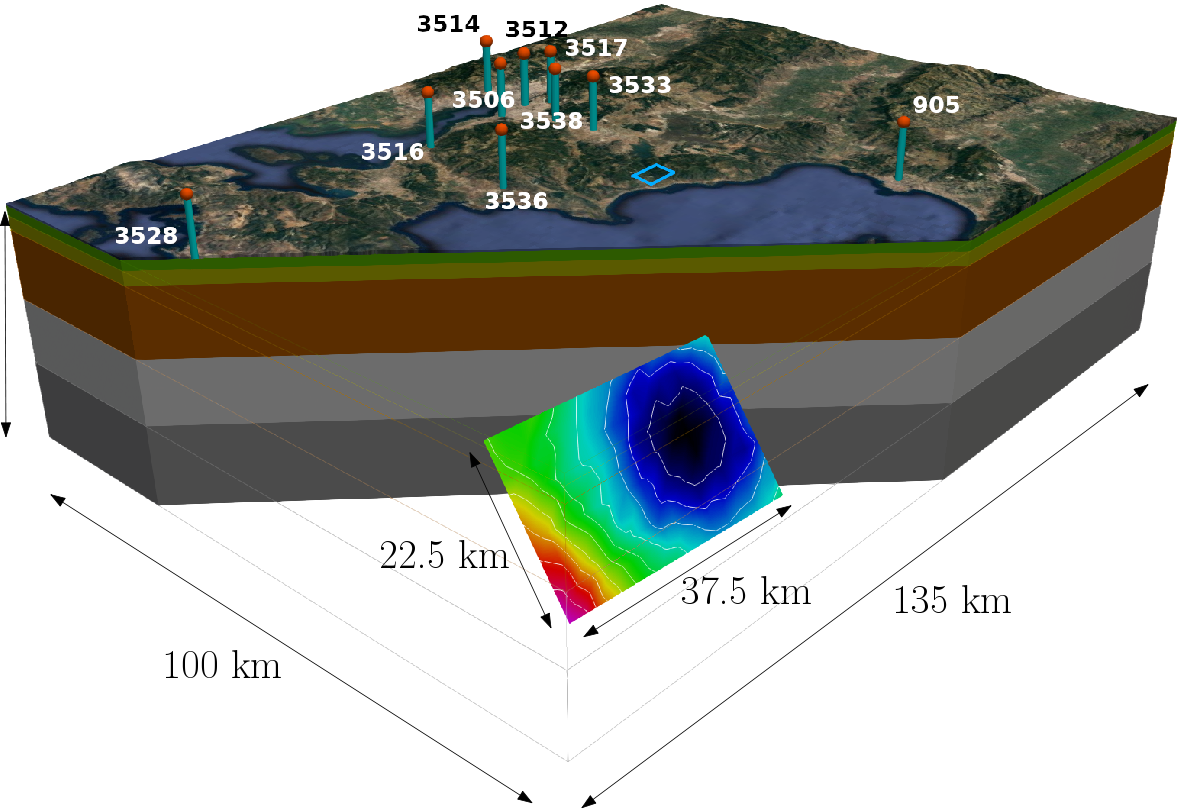}~~~\includegraphics[trim=0cm 0cm 0cm 0cm, clip, scale=0.175]{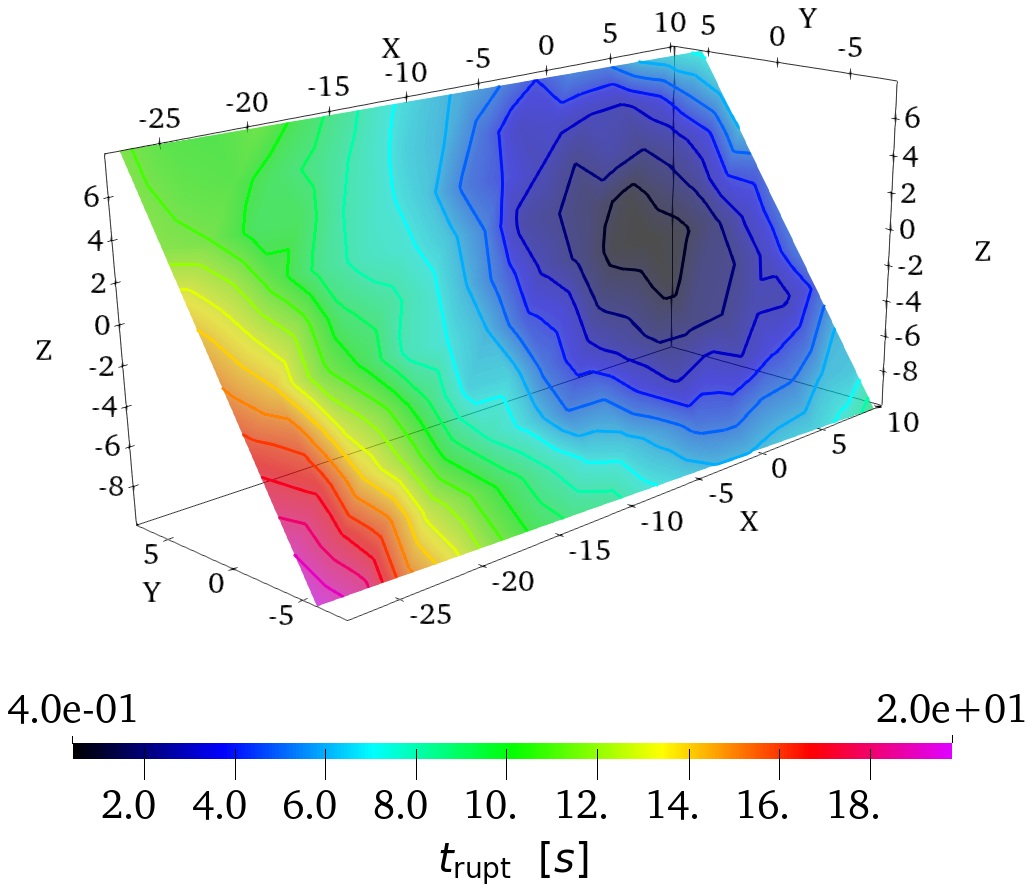}\\
		\vspace{-0mm}
		\hspace{2.2cm}\parbox{4.7cm}{\tiny Map data: Google, Imagery \textcopyright 2021 TerraMetrics, Map data \textcopyright 2021}
		\caption{\footnotesize Domain for the large scale simulation. Cross-section view with different material layers and the seismic fault plane showing iso-lines of rupture time (see also right figure for a zoom with time scale). The red pins mark the locations of the seismographic measurement stations with their AFAD ID numbers being printed next to them. The blue box is the dam vicinity region from Fig.\,\ref{Fig:Topography}.\label{Fig:LargeScaleDomain}}
	\end{center}
\end{figure}

For the aforementioned validation step, we first conduct a large scale simulation in the domain $\Omega_{\textup{large}}$, depicted in Fig.\,\ref{Fig:LargeScaleDomain}, without the fine scaled dam structure.\\

\paragraph{\bf Simulation parameters} The domain $\Omega_{\textup{large}}$ contains not only the area around the Tahtal{\i}-dam but also the seismic fault location where the considered event originated from (cf. Sec.\,\ref{sec:SeismicScenario}) as well as the locations of ten selected  stations of the AFAD network \cite{AFAD}, see Fig.\,\ref{Fig:LargeScaleDomain}.\kommentar{the two further real measurement stations \textit{AFAD\_905}, located at $37.86\,^{\circ}$N, $27.265\,^{\circ}$E \kommentar{{\color{gray}(21227.91863655 / -24511.06833 / 0)}} and \textit{Koeri\_GMLD}, located at $38.0766\,^{\circ}$N, $26.9157\,^{\circ}$E \kommentar{{\color{gray}(-9586.232823823 / -426.2472188 / 0)}} in addition to \textit{AFAD\_3536}.} The mesh used for this large scale simulation is coarser by a factor of around 13 (measured between two average elements with edge lengths 75\,m vs. 1000\,m), compared to the dam mesh from Fig.\,\ref{Fig:Blocking}, and ignores the dam structure as a detail in order to effectively span a volume of roughly $135 \times 100 \times 35$\,km. In order to properly describe the mechanical properties of the area along the Earth's crust, layered materials have been employed, which are assumed to be parallel to the $x$-$y$-plane. The adopted mechanical parameters as well as the layer depths used can be found in Tab.\,\ref{Tab:MaterialParamLargeSim}.

\begin{table}[h!]\renewcommand{\arraystretch}{1.4}
	\small
	\begin{tabular}{|c||c|c|c|c|c|c||c|}
		\cline{2-7}
		\multicolumn{1}{c||}{\textbf{Mat. par.}}&\multicolumn{6}{|c||}{\textbf{material block}} & \multicolumn{1}{c}{}\\ \hline
		\textbf{parameter} & \textbf{Layer 1} & \textbf{Layer 2} & \textbf{Layer 3} & \textbf{Layer 4} & \textbf{Layer 5} & \textbf{Layer 6} & \textbf{unit}\\ \hline\hline
		depth & $0-0.3$ & $0.3-1.7$ & $1.7-3.7$ & $3.7-13.7$ & $13.7-23.2$ & $23.2-35.0$ & km\\ \hline
		$\rho$ &  2355  &  2200   & 2300 &  2700 & 2900 & 3100 &  $\frac{\textup{kg}}{\textup{m}^3}$  \\
		 $v_p$  &  1695 & 2300  & 3200     & 6000& 6600&7200 & $\frac{\textup{m}}{\textup{s}}$  \\
		 $v_s$  & 1130   &  1600    &  3400     & 3700 & 3700& 4000& $\frac{\textup{m}}{\textup{s}}$  \\
		$Q_s$ & 113 & 160 & 340 & 370 & 370 & 400 & $-$  \\ \hline
	\end{tabular}	
	\caption{\footnotesize \textbf{Layers and material parameters} used for the large-scale simulations. Data derived and adapted from the original data-set taken from \cite{USGS}. The quality factor $Q_s =  \pi f_0/\zeta$, where $f_0$ is a frequency reference value here chosen equal to $f_0=1$\,Hz. }\label{Tab:MaterialParamLargeSim}
\end{table}

The final large-scale mesh contains $N_{\textup{el,large}}=536.105$ elements with polynomial degree of $p=3$ for the ansatz functions. Time discretization is done with $N_{T,\textup{large}}=40.000$ timesteps of size $\Delta t_{\textup{large}}=10^{-3}~s$.\\

\begin{figure}[h!]
	\begin{center}
		\includegraphics[trim=0cm 0cm 0cm 0cm, clip, scale=0.2]{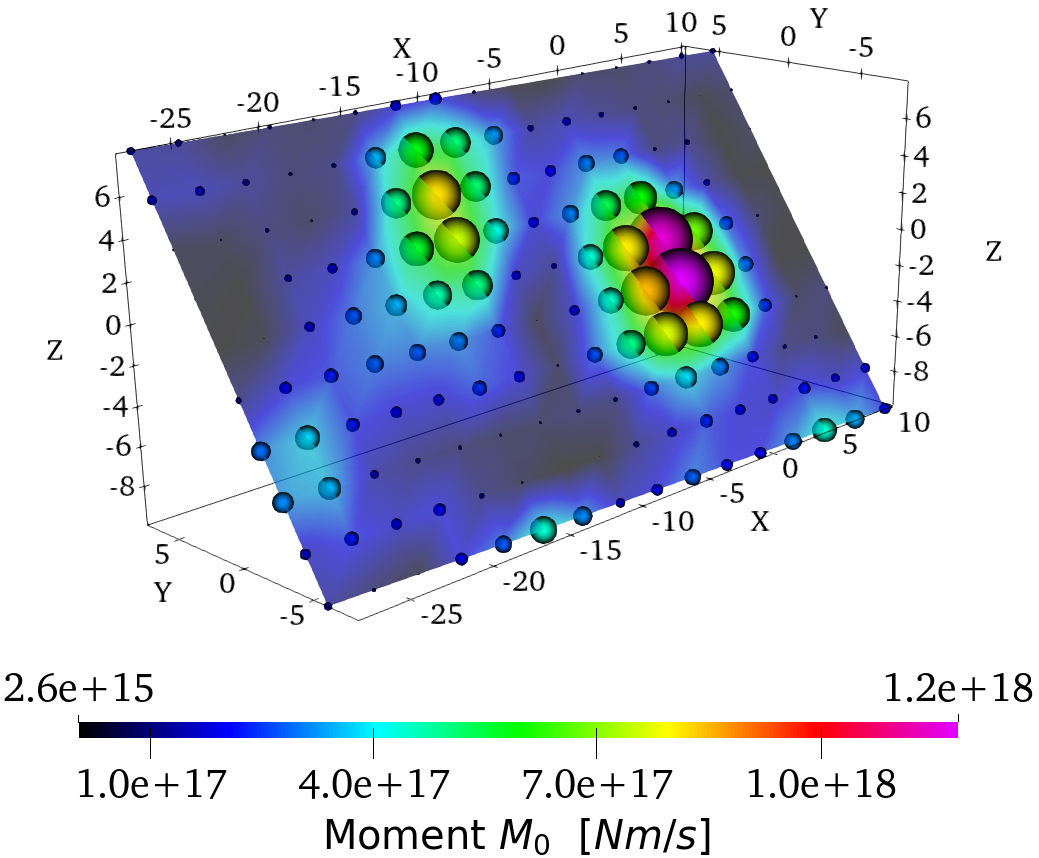}\includegraphics[trim=0cm 0cm 0cm 0cm, clip, scale=0.2]{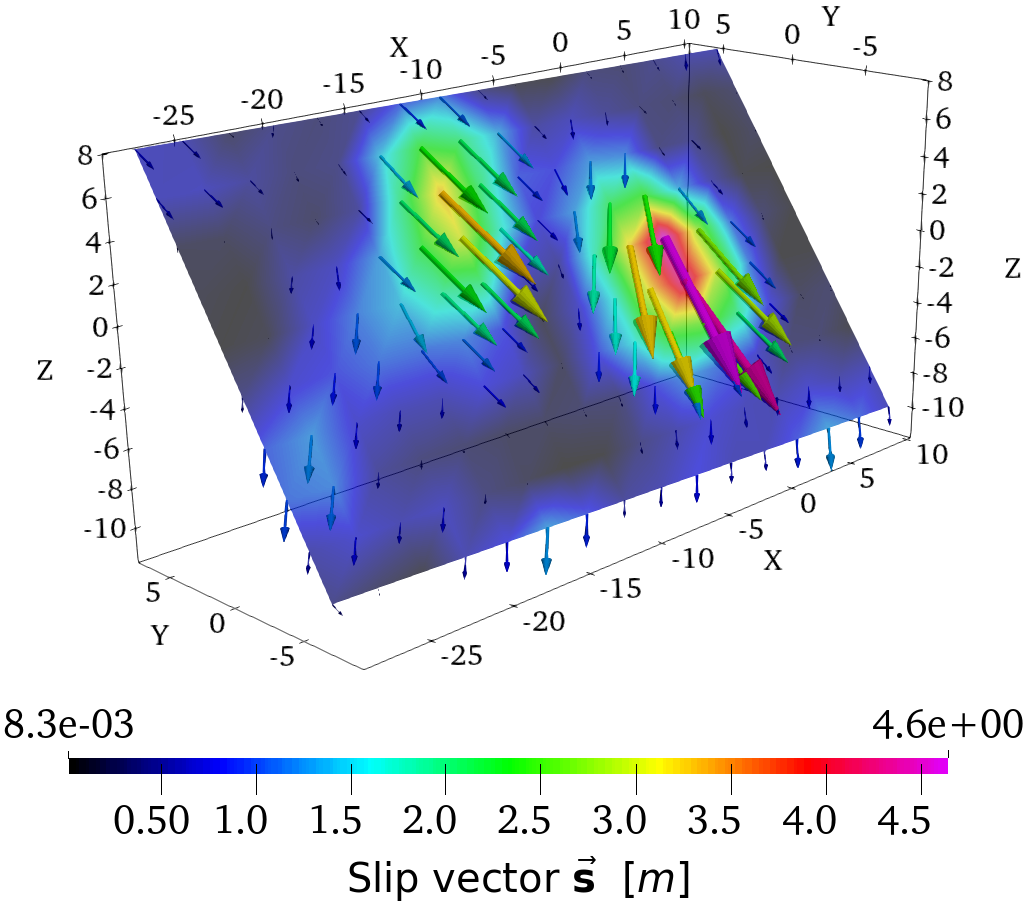}\\
		\caption{\footnotesize Seismic fault plane with coordinates relative to the hypocenter measured in km. \textbf{(left)} Seismic moment magnitude $M_0$, \textbf{(right)} slip vector field $\vecc{s}$ (vectors scaled by a factor of 2000), rupture time $t_{\textup{rupt}}$ can be found in Fig.\,\ref{Fig:LargeScaleDomain}.\label{Fig:SourceData}}
	\end{center}
\end{figure}

\paragraph{\bf Seismic fault data and kinematic source mechanism} The source mechanism is described by means of a set of double-couple moment-tensors  $\tens{M}_i(t), i=1,2,\dots,160$ distributed along the fault plane, each with its own set of source parameters being slip-vector $\vecc{s}_i$, rupture- and rise-times $t_{\textup{rup},i}$ and $t_{\textup{rise},i}$ as well as the released moment magnitude $M_{0,i}$. Data have been obtained by \cite{USGS} on a grid of $10\times 16$ points $p_i$ across the seismic fault plane depicted in Fig.\,\ref{Fig:SourceData}.
The moment tensors $\tens{M}_i$ are then associated with the numerical quadrature node closest to the data point $p_i$ and can be computed as:
\begin{equation*}
	\tens{M}_i(t) = M_{0,i}\cdot m_{i}\left(\frac{t- t_{\textup{rupt},i}}{t_{\textup{rise},i}}\right)\cdot \left[\left(\vecc{s}_i\otimes \vecc{n}\right)+\left(\vecc{s}_i\otimes \vecc{n}\right)^{\top}\right]
\end{equation*}
where $\vecc{n}$ is the fault plane normal computed from the provided data-set, and $m_{i}(\hat{t})$ is a normalized moment-function monotonically increasing from 0 to 1 that models the moment-release over time at the point $p_i$ respecting the available data of rupture- and rise-time, hence the individual distance from the hypocenter. We refer the reader to \cite{archuleta2013} for the precise definition of the moment-rate functions $\dot{m}_{i}$ from which $m_i$ are computed.\\

\begin{figure}[h!]
	\begin{center}
		\includegraphics[trim=0cm 0cm 0cm 0cm, clip, scale=0.25]{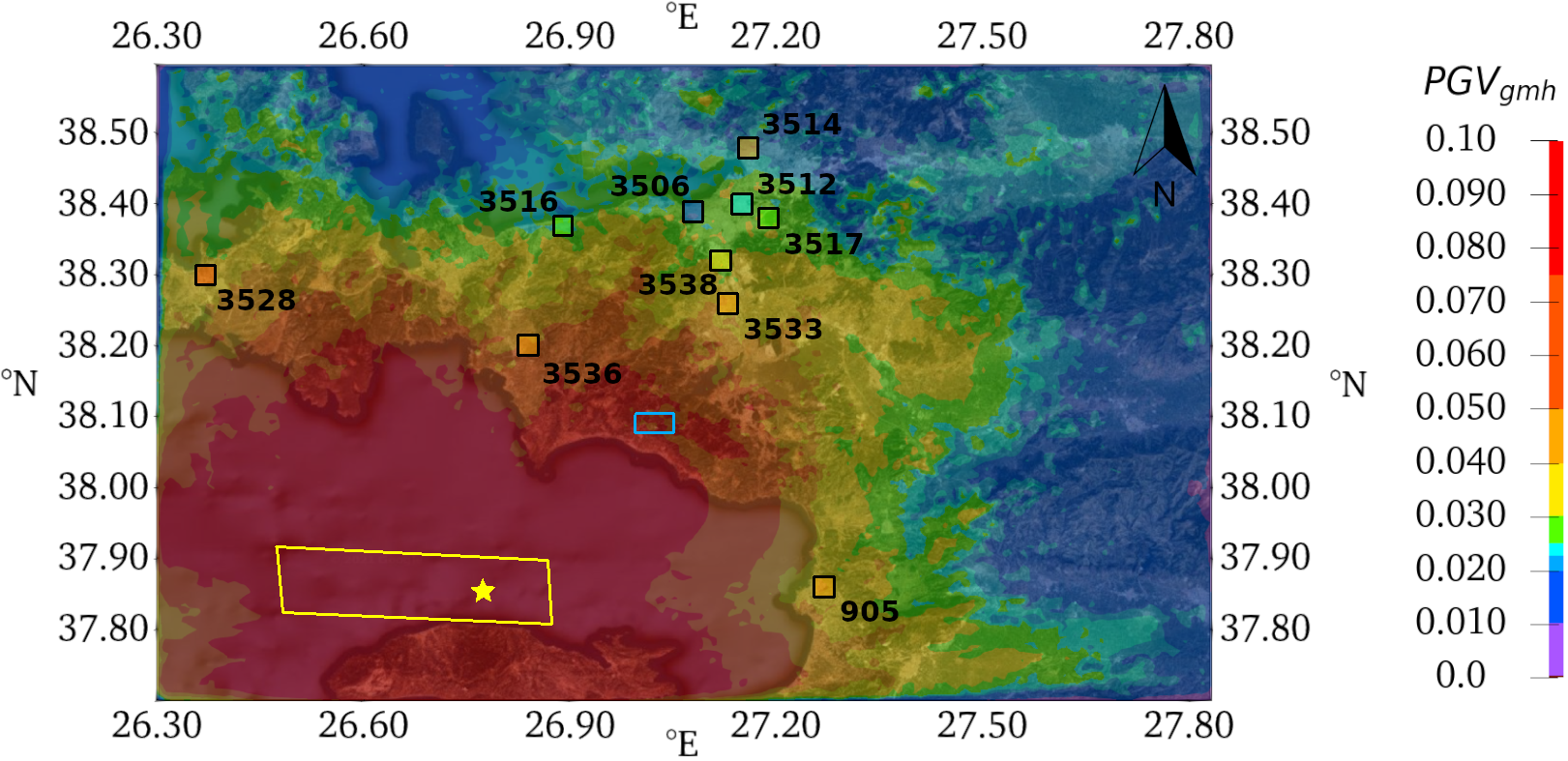}\\
		\vspace{-0mm}
		\hspace{8.2cm}\parbox{4.7cm}{\tiny Map data: Google, Imagery \textcopyright 2021 TerraMetrics, Map data \textcopyright 2021}
		\caption{\footnotesize PGV-map for the large scale simulation. Color scale is cut-off at $0.1\,\frac{\textup{m}}{\textup{s}}$. AFAD stations are marked by their ID number, peak ground velocity values of measurements are color coded. The blue rectangle shows the dam vicinity region of Fig.\,\ref{Fig:Domain}, the hypocenter and the fault are marked in yellow. \textit{(original map image overlayed with simulation colormap)}\label{Fig:LargeUVOutput}}
	\end{center}
\end{figure}

\begin{figure}[h!]
	\begin{center}
		\includegraphics[trim=0cm 0cm 0cm 0cm, clip, scale=0.4]{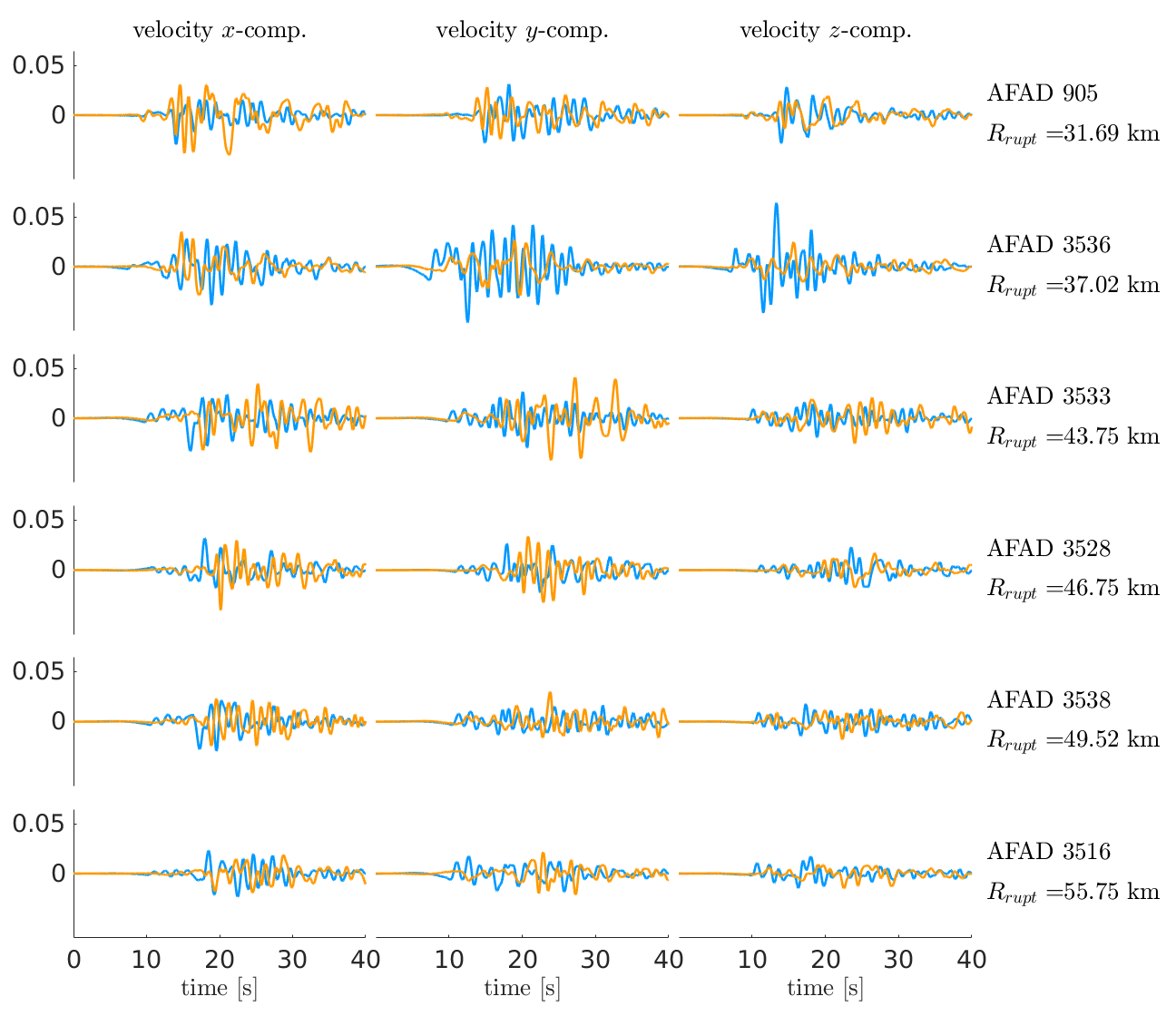}
		\caption{\footnotesize \footnotesize\textbf{(rows)} Comparison of numerically evaluated velocity seismograms ({\color{MATblue}blue}) and records ({\color{orange}orange}) for the six AFAD-stations closest to the hypocenter. Stations are sorted by increasing $R_{\textup{rupt}}$ (see label on the right, also for station ID). \textbf{(columns)} First column $x$-component, second $y$-component, third $z$-component. All data filtered to $[0.1,1]\,\textup{Hz}$ with $2^{\textup{nd}}$ order Butterworth filter. \label{Fig:LargeComparison}}
	\end{center}
\end{figure}

\begin{figure}[h!]
	\begin{center}
		\includegraphics[trim=0cm 0cm 0cm 0cm, clip, scale=0.4]{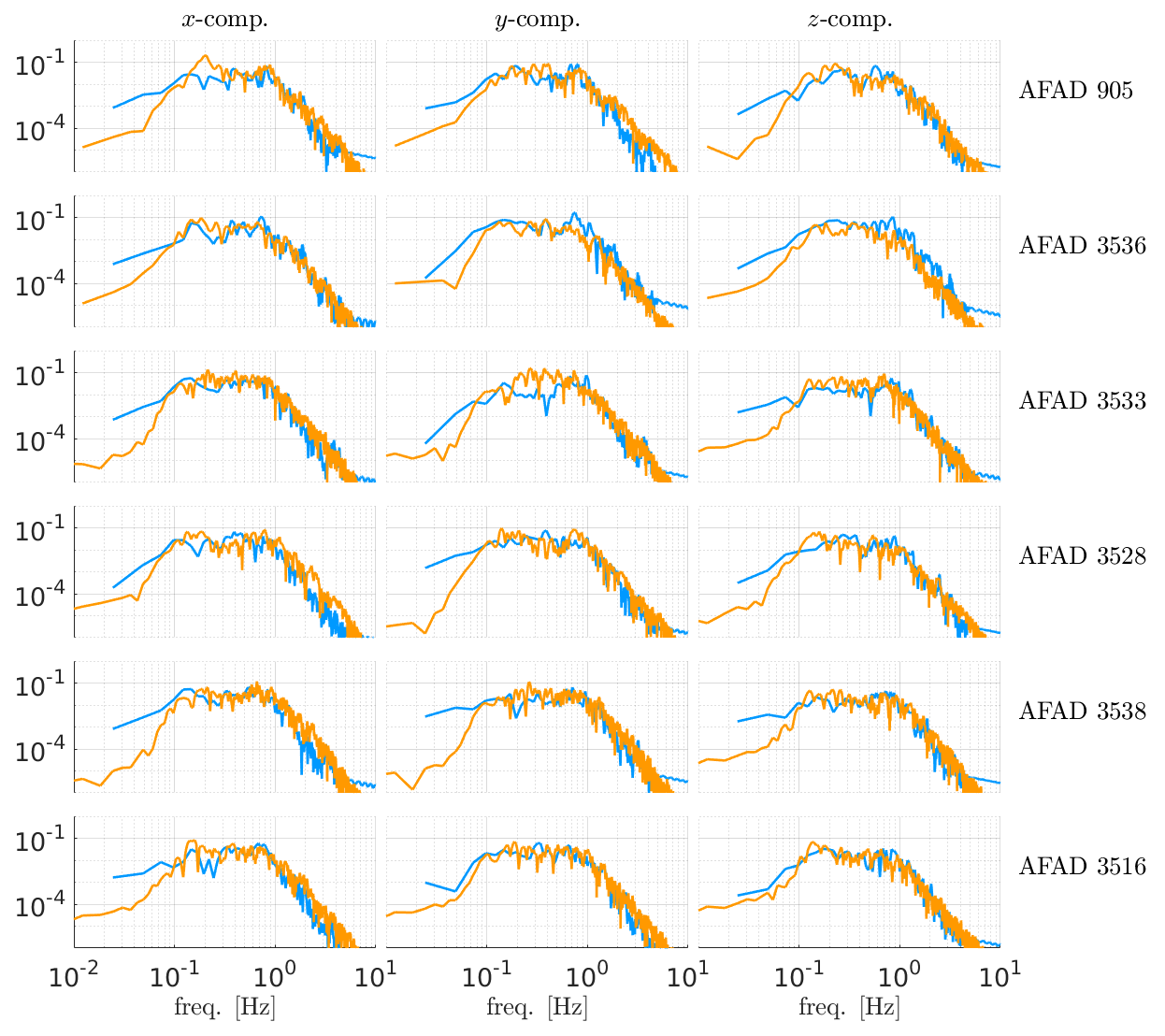}
		\caption{\footnotesize Amplitude power spectra of the velocity signals in Fig.\,\ref{Fig:LargeComparison}. Again ({\color{MATblue}blue}) are the numerical data, ({\color{orange}orange}) the filtered measurement data. \label{Fig:LargeFourierComparison}}
	\end{center}
\end{figure}

\paragraph{\bf Numerical results and validation} In order to assess the maximal (over time) displacement, resp. velocity that is attained at each point $\vecc{p}$ on the computational domain's surface, we introduce the so called \textit{\textbf{g}eometric \textbf{m}ean \textbf{h}orizontal peak ground displacement} $\pgu$ and \textit{velocity} $\pgv$ at point $\vecc{p}\in\Omega$ \cite{BooreBommer2005} as quantities of interest:
{\small \begin{align*}
\pgu(\vecc{p}):=\sqrt{\sup_{t\in(0,T)}\vecc{u}_x(t,\vecc{p})\cdot \sup_{t\in(0,T)}\vecc{u}_y(t,\vecc{p})},\quad \pgv(\vecc{p}):=\sqrt{\sup_{t\in(0,T)}\vecc{v}_x(t,\vecc{p})\cdot \sup_{t\in(0,T)}\vecc{v}_y(t,\vecc{p})}
\end{align*}}
Fig.\,\ref{Fig:LargeUVOutput} then shows the PGV-map of the large scale simulation. The location of the dam, fault and hypocenter as well as the AFAD stations available are reported (see also Fig.\,\ref{Fig:LargeScaleDomain}). The latters are color coded, according to the values retrieved from the observed seismograms. The filtered time-history of simulated and measured seismograms are also listed and compared in the time domain in Fig.\,\ref{Fig:LargeComparison} and in the frequency domain in Fig.\,\ref{Fig:LargeFourierComparison}. Both for the subset of the six stations closest to the hypocenter.

\begin{figure}[h!]
	\begin{center}
		\includegraphics[trim=0cm 0cm 0cm 0cm, clip, scale=0.46]{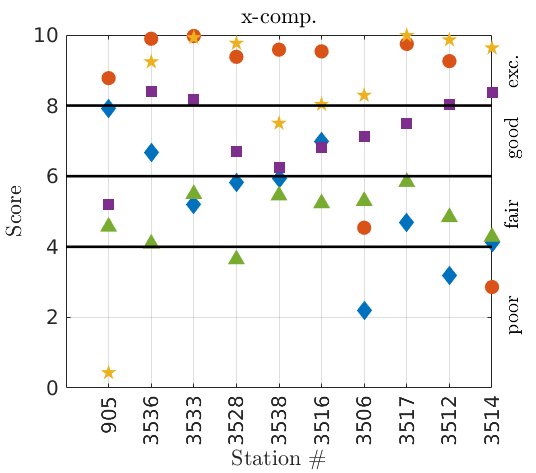}\hspace{5mm}\includegraphics[trim=0cm 0cm 0cm 0cm, clip, scale=0.46]{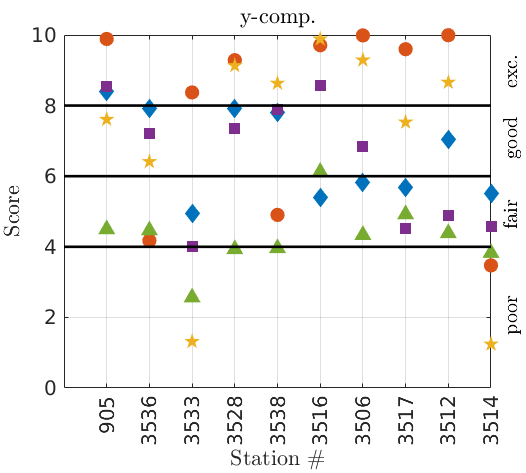}\includegraphics[trim=0cm 0cm 0cm 0cm, clip, scale=0.46]{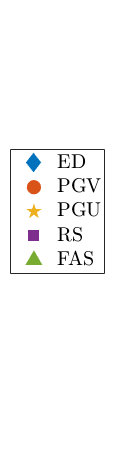}\\[5mm]
		\includegraphics[trim=0cm 0cm 0cm 0cm, clip, scale=0.46]{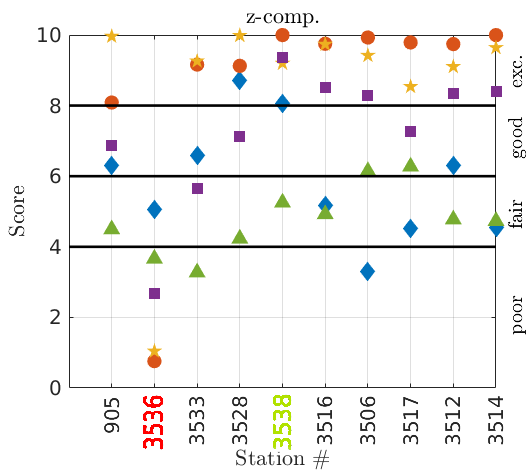}\hspace{5mm}\includegraphics[trim=0cm 0cm 0cm 0cm, clip, scale=0.2]{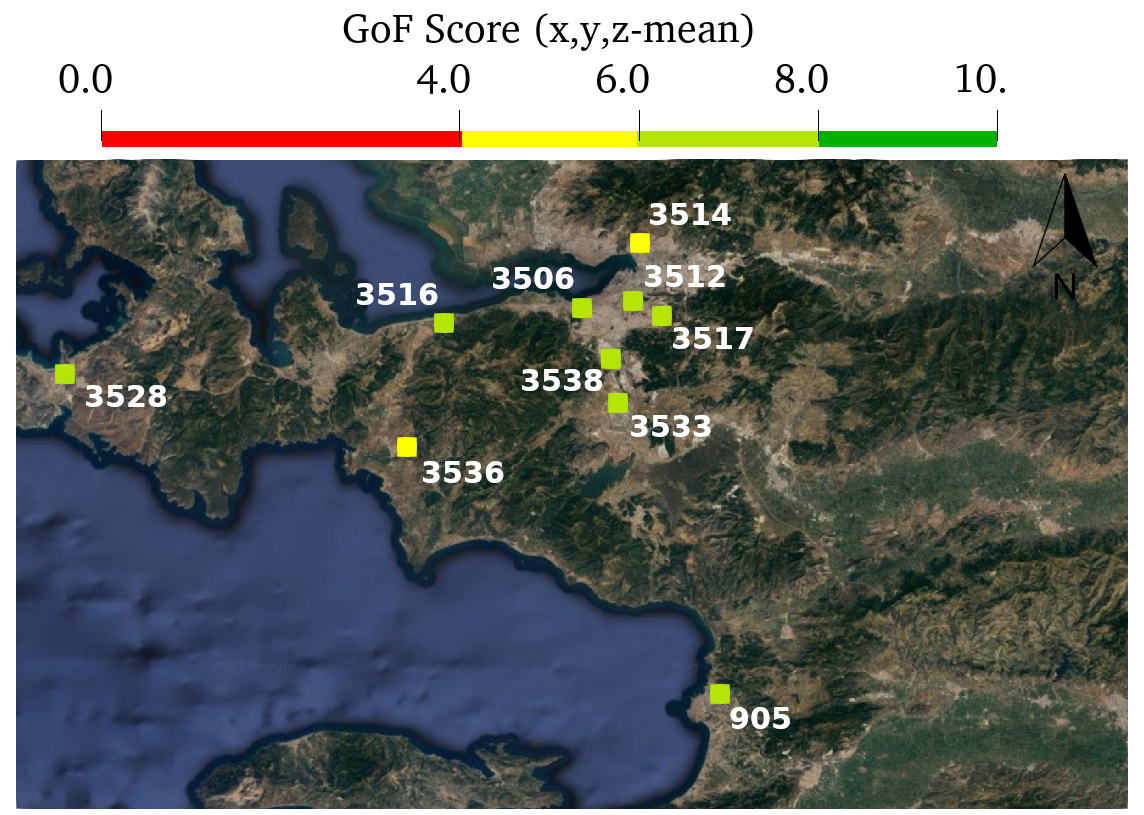}\\
		\vspace{-0mm}
		\hspace{8.2cm}\parbox{4.7cm}{\tiny Map data: Google, Imagery \textcopyright 2021 TerraMetrics, Map data \textcopyright 2021}
		\caption{\footnotesize Evaluation of Anderson GoF-criteria to assess the similarity of the seismographs in Fig.\,\ref{Fig:LargeComparison}. \textbf{(top left, right, bottom left)} Energy duration (ED), peak ground velocity (PGV), peak ground displacement (PGU), response spectrum (RS) and Fourier amplitude spectrum (FAS) evaluated for $x,y$ and $z$-components at the presented stations. Extreme values marked in red and green: $z$-comp. of \#\,3536 at $\approx 2.637$ (poor) to $z$-comp. of \#\,3538 at $\approx 8.372$ (excellent). \textbf{(bottom right)} Station-wise mean of $x,y$ and $z$ component's presented Anderson criteria.\label{Fig:Anderson}}
	\end{center}
\end{figure}

Given the relatively simplistic model adopted, Fig.\,\ref{Fig:LargeComparison} shows a satisfactory agreement between the observed and simulated velocities, in terms of arrival time, duration of the signal, phase and amplitude of the waves. The goodness of these results is also confirmed by the PGV map, showed in Fig.\,\ref{Fig:LargeUVOutput}. Referring to the comparison in terms of Fourier spectra, cf. Fig. \ref{Fig:LargeFourierComparison}, it turns out that, in general, there is a satisfactory agreement between simulated and recorded amplitudes for frequencies up to about 1\,Hz, although synthetic tends sometimes to underestimate the observed amplitudes.\\
\indent To give a quantitative measure of the overall performance of the numerical simulation we adopt the Goodness of Fit (GoF) criteria proposed by  \cite{Anderson2004}, being this latter widely followed and recognized for this kind of evaluation. For the frequency band of interest (i.e. $0.1-1$\,Hz), a GoF score from 0 to 10 ($<4$, poor; $4-6$, fair; $6-8$, good; $\geq 8$, excellent) is estimated on five metrics of interest for engineering purposes, namely: energy duration (ED), Peak Ground Velocity (PGV), Peak Ground Displacement (PGU), Response Spectral (RS) acceleration and Fourier Amplitude Spectrum (FAS). Note that FAS and RS criteria are evaluated considering only the frequencies and structural periods within the range $0.1-1$\,Hz of validity of the numerical simulations. The GoF scores, computed for each criterion and for the three components of motion, are shown in Fig. \ref{Fig:Anderson} for the whole set of ten recording stations considered in Fig.\,\ref{Fig:LargeUVOutput}. These results confirm that with few exceptions, the numerical model provides predictions that are in overall good agreement (from fair to excellent) with the records. By taking inspiration from \cite{Olsen2010} we not only compute the average GoF value for each ground motion component but also the average between the components in order to summarize the results into a single final score for each station. Hence, in Fig. \ref{Fig:Anderson}, we additionally present that score as an overview of the misfit between records and simulated results at territorial scale. We found that for the majority of the stations the agreement between simulations and observations is from fair to good. These results are aligned with those obtain for other, different earthquake scenarios as, e.g., \cite{Paolucci2015, Guidotti2011}.

\subsection{Plane-wave excitation using nearby seismogram}\label{subsec:PlaneWave1}
\paragraph{\bf Measurement data and source mechanism} For our second numerical simulation, we use input data for the ground motion from recorded measurements. Namely the station AFAD \# 3536, located at $38.1968\,^{\circ}$N, $26.8384\,^{\circ}$E \kommentar{{\color{gray}(-16405.39660248744 / 12939.38296387250 / 0)}}, (cf. Figs.\,\ref{Fig:LargeScaleDomain},\,\ref{Fig:LargeUVOutput}) and positioned approximately 30 km from the Tahtal{\i}-dam. In Fig.\,\ref{Fig:LargeComparison}, second row, the $x,y$ and $z$ components $v_{x,y,z}^{\textup{ref}}(t)$ (orange lines) of the ground motion velocity $\vecc{v}^{\textup{ref}}(t)$ recorded at AFAD \# 3536 during the seismic event \cite{AFAD} can be seen. 
They were deconvolved to a certain depth, since it was observed on the free surface, and subsequently it was adopted as input for the plane wave excitation, following \cite{FaccioliQuarteroni}. Therein
the equivalent body force 
\begin{equation*}
	\vecc{f}=2\rho_{\textup{e}}v_p\delta(z-z_0)\vecc{v}^{\textup{ref}}(t)
\end{equation*}
is applied on a horizontal plane located at $z=z_0$. The seismic wave generated by it will rise and yields an approximation to the actual ground motion of the event. The realization of this approach is based on an additional layer of same material properties below the bottom surface of the mesh in Fig.\,\ref{Fig:Blocking}, right, which has a thickness of only one element and $z_0$ being its mean depth, in which the body force is applied. This approach is often used in engineering analysis aiming at simulating the so-called ``dynamic soil structure interaction'' problem (DSSI). In fact, due to scarcity of numerical code capable of taking into account the entire problem (from source-to-site) and/or having only limited computational resources available, the state-of-the-art engineering approach considers, typically, only the region in the immediate proximity to the dam and assumes that at this scale the excitation can be properly approximated by a plane wave \cite{Sunbul2017, Bayraktar2009b}. Obviously this approach presents a series of limitations that are even more relevant when, as often happens operationally, the dam is studied in 2D \cite{Bayraktar2009a}.\\

\paragraph{\bf Simulation parameters} Tab.\,\ref{Tab:MaterialParam} shows the material parameters used for the different zones of the numerical simulations. Since "in-situ" measured values were not available it was decided to adopt reliable literature values.

The numerical simulation was conducted with $N_T=10^6$ timesteps of size $\Delta t=5\cdot 10^{-5}$ on a spacial grid using $N_{\textup{el}}=18.948$ elements with a polynomial degree of $p=2$ for the ansatz functions. The DG penalty-parameter was chosen as $\beta=250$.\\
\begin{table}[h!]\renewcommand{\arraystretch}{1.4}
	\normalsize
	\begin{tabular}{|c|c||c|c|c||c|}
		\cline{3-5}
		\multicolumn{2}{c||}{\textbf{Material parameters}}&\multicolumn{3}{|c||}{\textbf{material block}} & \multicolumn{1}{c}{}\\ \hline
		\textbf{type} & \textbf{parameter} & \textbf{Ground} & \textbf{Dam} & \textbf{Water} & \textbf{unit}\\ \hline\hline
		all    & $\rho$ & 2355   & 2000    & 998.23 & $\frac{\textup{kg}}{\textup{m}^3}$\\ \hline
		& $v_p$  & 1695 & 525  &   /    & $\frac{\textup{m}}{\textup{s}}$\\
		elastic& $v_s$  & 1130   & 350     &   /    & $\frac{\textup{m}}{\textup{s}}$\\
		&$Q_s$ & 500    &      500    &   /    & - \\ \hline
		\raisebox{-5mm}{acoustic}& $c$     &   /    &   /   & 1500 & $\frac{\textup{m}}{\textup{s}}$\\
		& $b$     &   /    &   /   & $6\cdot 10^{-9}$ &$\frac{\textup{m}^2}{\textup{s}}$\\ \hline
	\end{tabular}~\\[5mm]
	\caption{\footnotesize \textbf{Material parameters} used for the numerical simulations. Ground consists of Rock/Limestone and the Dam is mostly modelled as Gravel/Sand material \label{Tab:MaterialParam}. The quality factor $Q_s=\pi f_0/\zeta $, where $f_0$ is a frequency reference value here chosen equal to $f_0= 1$\,Hz.}
\end{table}

\begin{figure}[h!]
	\begin{center}
	    \textbf{Plane wave:}\parbox{12cm}{~}\\
		\includegraphics[trim=0cm 0cm 0cm 0cm, clip, scale=0.14]{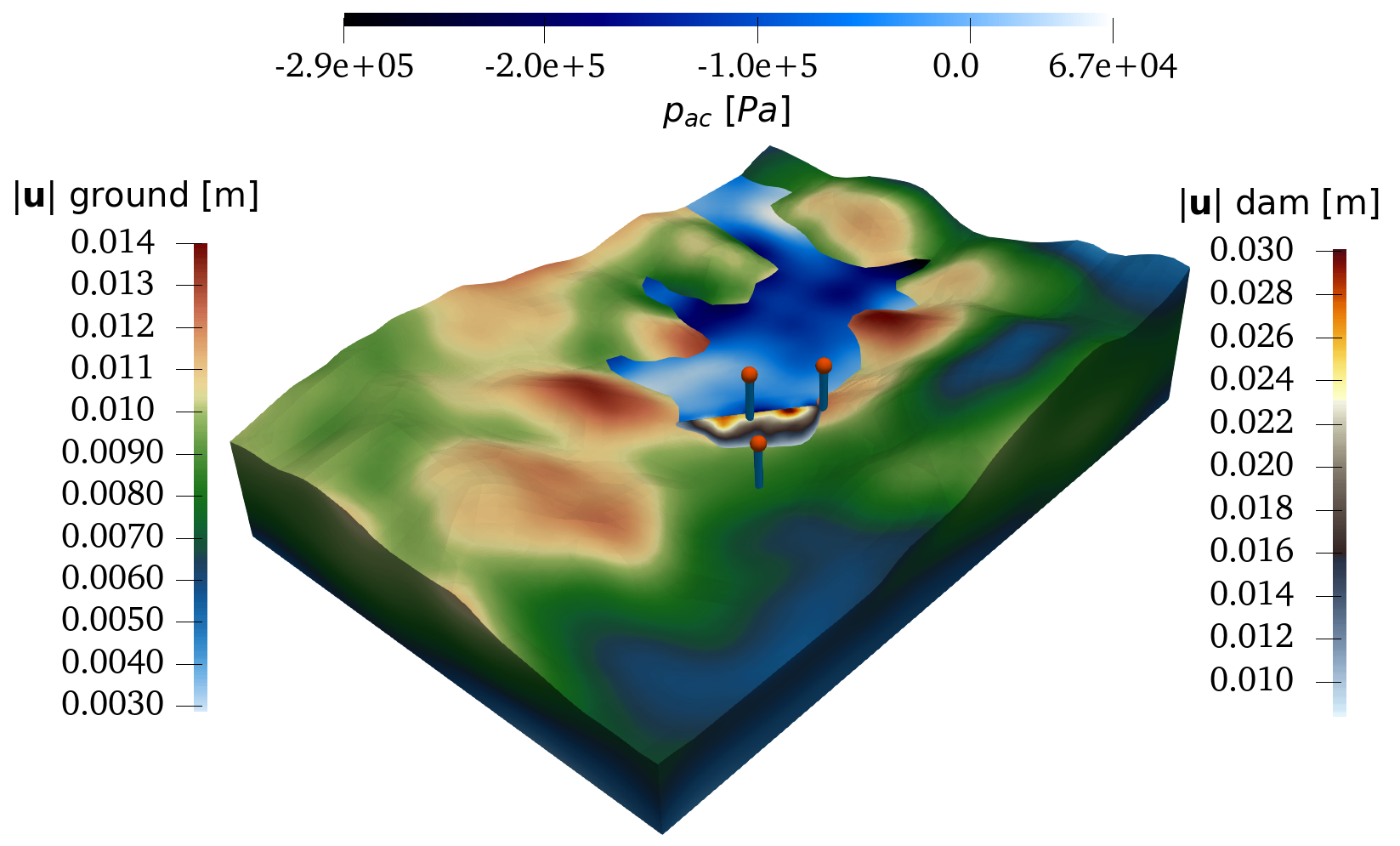}\\
		\textbf{Source-to-site:}\parbox{12cm}{~}\\
		\includegraphics[trim=0cm 0cm 0cm 0cm, clip, scale=0.2]{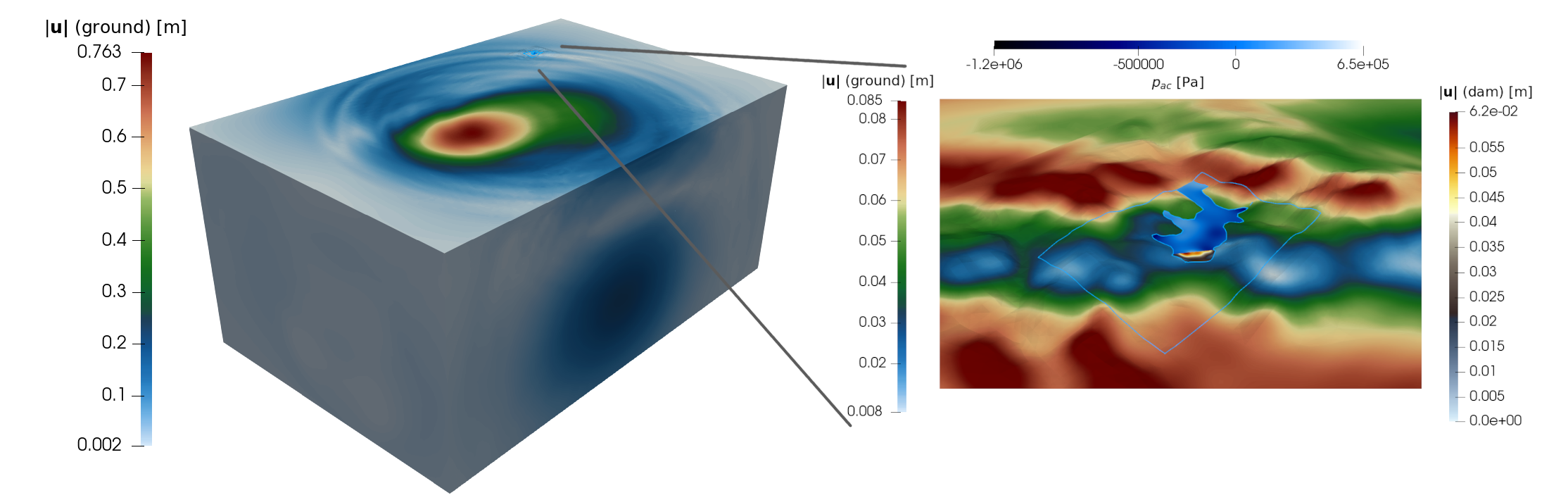}\\
		\caption{\footnotesize Time snapshot of displacement-field in \textbf{(top)} plane wave simulation, \textbf{(bottom)} source-to-site simulation where in the zoomed in picture the blue frame of the dam-vicinity region $\Omega$ is visible. The pin-needles (top picture) show the locations of the synthetic seismographs evaluated in Fig.\,\ref{Fig:PlaneWave2vsS2S}.\label{Fig:UVA_Field}}
	\end{center}
\end{figure}

\paragraph{\bf Numerical results} Fig.\,\ref{Fig:UVA_Field} (top) shows a snapshot of the simulated displacement-field $\vecc{u}$ in the elastic, and acoustic pressure field $p_{\textup{ac}}$ in the acoustic domain. The bottom row shows a corresponding snapshot of the upcoming source-to-site simulation from Sec.\,\ref{Sec:S2S} which highlights the different scales of magnitude of the simulations. In addition to the time-snapshot picture, we can also employ the time-history of the simulation at certain points of interest (cf. Fig.\,\ref{Fig:PlaneWave2vsS2S}) and again compute peak-ground maps of the conducted simulation, which will also be depicted later in Fig.\,\ref{Fig:PlaneWave2vsS2SPG} for a comparative discussion with the results of the full source-to-site simulation.
Fig.\,\ref{Fig:PlaneWave2vsS2S} contains the numerical seismograms of the components of $\vecc{v}$ (oriented orthogonal, parallel and vertical to the dam, compare Fig.\,\ref{Fig:DamDisplacementSketch} for orientation) of the simulation at three distinct locations around the dam. One of them is located on the ground in front of the dam\kommentar{{\color{gray}(1239.583 / 673.79655 / 13.180415)}}, the other directly on top of the dam\kommentar{{\color{gray}(1478.465 / 928.53832 / 73.5)}} and the third at the right dam ambutment. The exact locations of these synthetical seismographs are depicted in Fig.\,\ref{Fig:UVA_Field} via the red pin-needles. In order to compare the signals with similar frequency content, the signals are again frequency filtered based on a second order Butterworth filter with a frequency-band of $[0.1, 1]$\,Hz, eliminating high frequency (numerical) artifacts and measurement oscillations.


\subsection{Full source-to-site simulation}\label{Sec:S2S}
\begin{figure}[h!]
	\begin{center}
		\hspace*{-5mm}\includegraphics[trim=0cm 0cm 7.5cm 0cm, clip, scale=0.12]{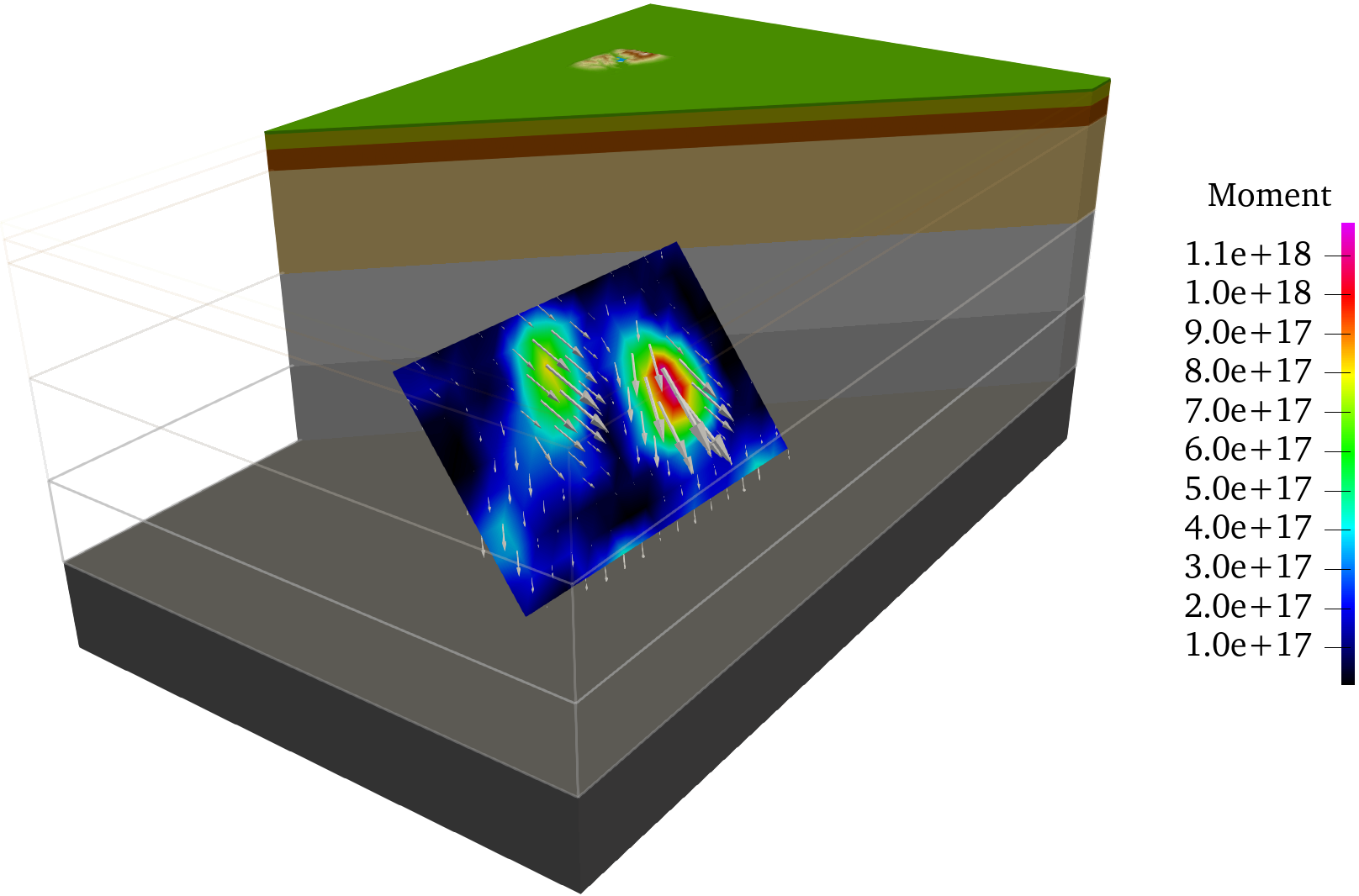}\hspace*{1cm}\raisebox{7mm}{\includegraphics[trim=0cm 0cm 0cm 0cm, clip, scale=0.09]{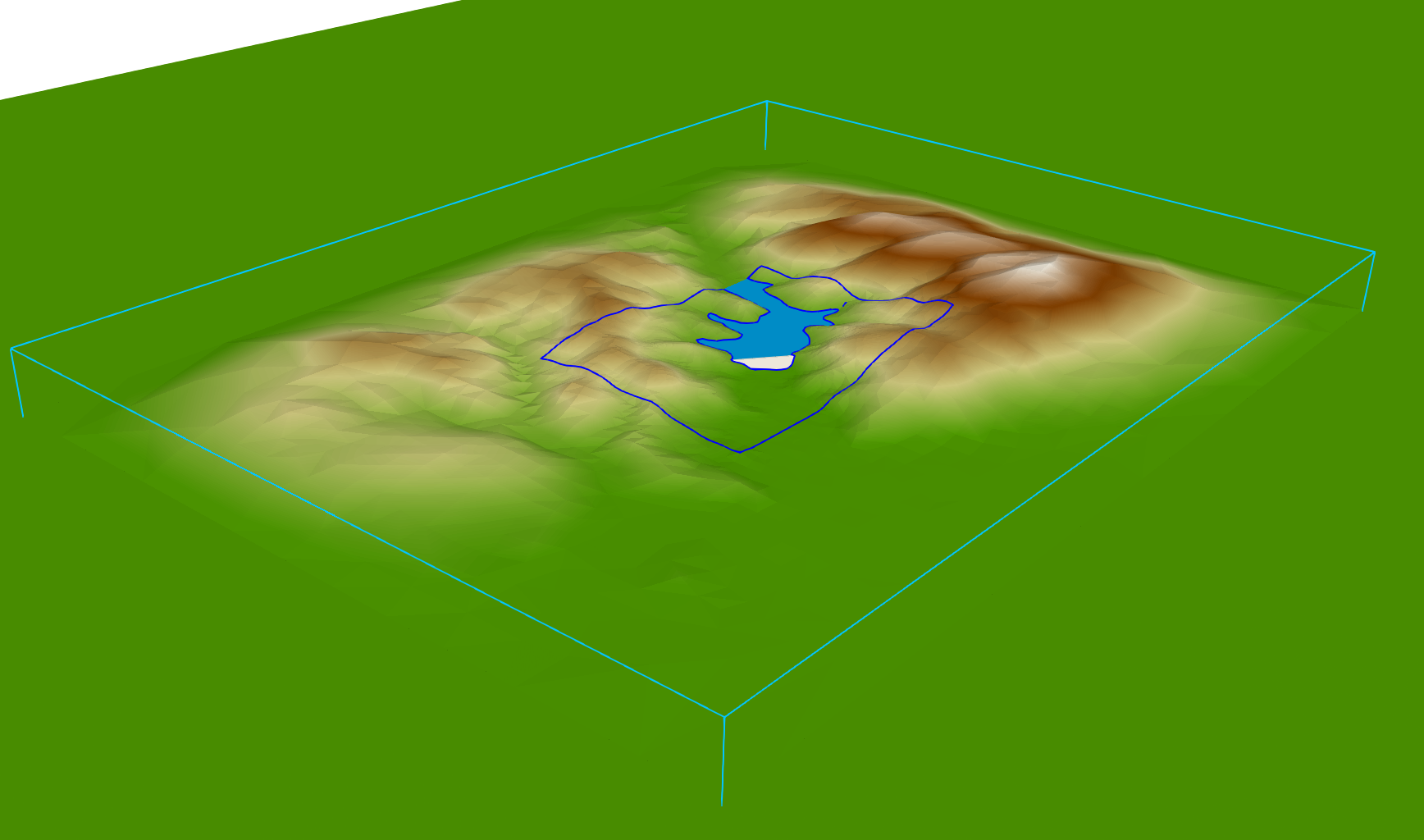}}\\
		\includegraphics[trim=0cm 0cm 0cm 0cm, clip, scale=0.16]{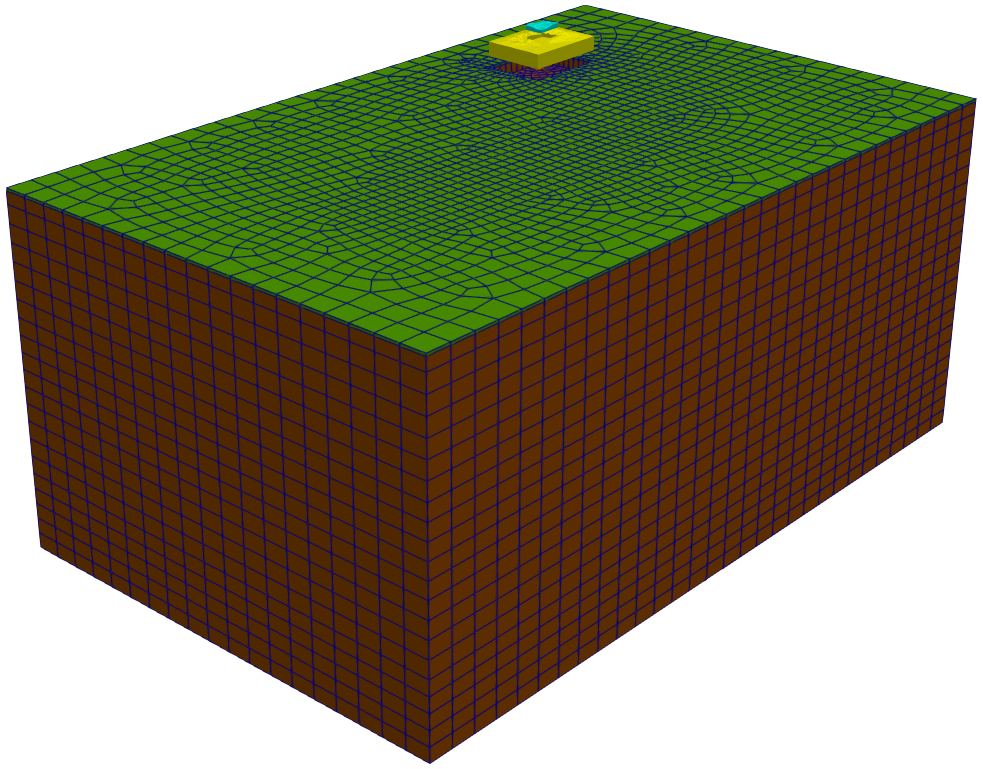}\hspace*{1cm}\raisebox{7mm}{\includegraphics[trim=0cm 0cm 0cm 0cm, clip, scale=0.12]{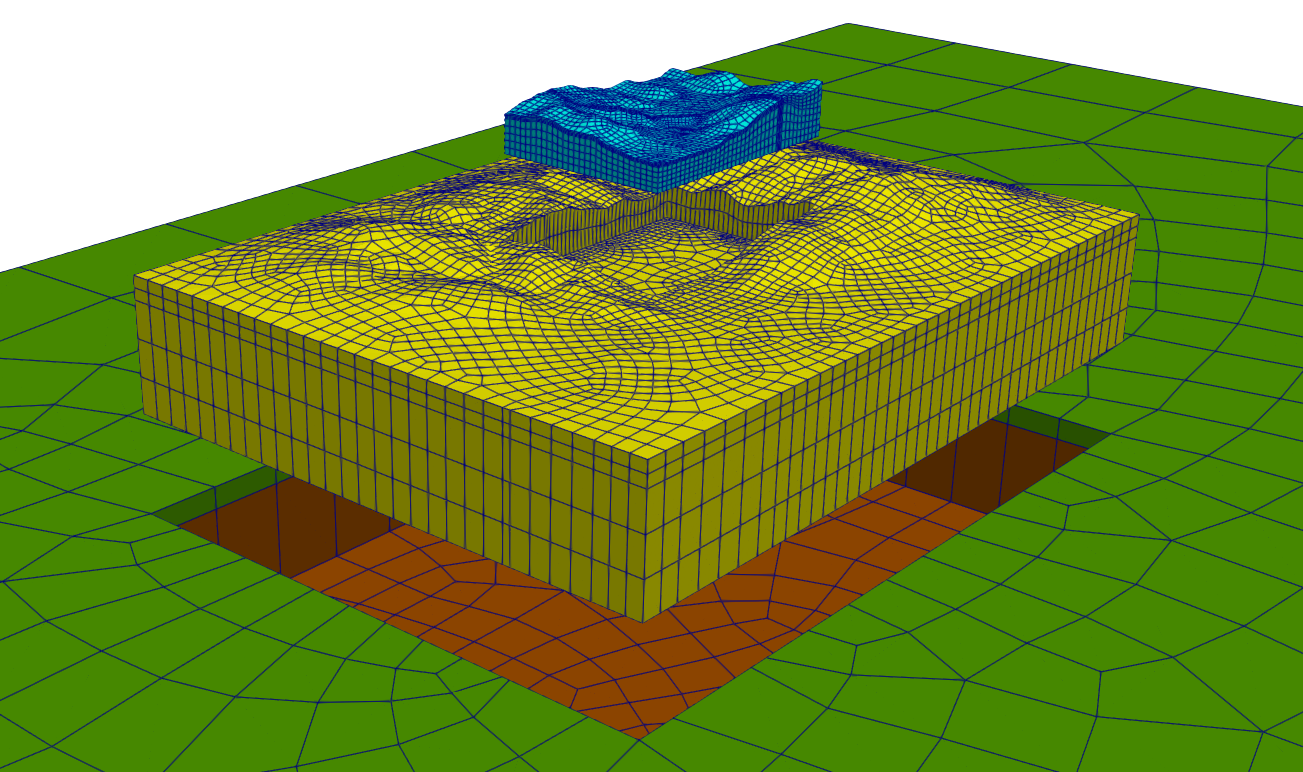}}
		\caption{\footnotesize\textbf{(top left)} Full source-to-site domain in cross-section view with different material layers as well as the seismic fault plane and slip distribution (cf. Fig. \ref{Fig:SourceData}) visible. \textbf{(top right)} Zoom into the region of the dam. The dark blue lines mark $\Omega$ as a sub-region with resolved topography (cf. Fig. \ref{Fig:Domain}), the light blue frame marks the transition zone. \textbf{(bottom left)} Hexahedral mesh of the complete source-to-site domain $\Omega_{\textup{s2s}}$. The brown surfaces carry absorbing boundary conditions. \textbf{(bottom right)} Zoom into the region of the dam. The transition zone (yellow) and the core blocks (blue) with their individual meshes are elevated here in order to see the (non matching) DG interfaces between each of them.\label{Fig:S2SDomain}}
	\end{center}
\end{figure}

The final simulation setup consists of a full source-to-site simulation using a domain $\Omega_{\textup{s2s}}$ containing the Tahtal{\i}-dam area as well as the seismic fault just as $\Omega_{\textup{large}}$. By respecting the dam-structure, the topography in its vicinity and the presence of water behind the dam, $\Omega$ becomes a subset of $\Omega_{\textup{s2s}}$ and hence the original mesh (cf. Fig.\,\ref{Fig:Blocking}) a submesh as well. The resulting grid contains multiple length-scales in element size from the larger elements having an edge length in the order of a kilometer to the smaller ones measuring only few meters. The geometry is created by embedding the domain $\Omega$ from Fig.\,\ref{Fig:Domain} on top of a rectangular block of an approximate size of $88\times55\times40\,$km. In order to have a smooth top-surface, a transition zone was used, in which the detailed topography of the dam-area gets flattened down to the level $z=0$ such that it is accurately resolved only in the proximity of the dam. Fig.\,\ref{Fig:S2SDomain} shows the resulting geometry, which was then meshed with a locally refined mesh consisting out of around 100.000 elements in total. The mesh is also depicted in Fig.\,\ref{Fig:S2SDomain} and takes advantage of the DG coupling approach by individually meshing the surrounding, transition and core blocks in a not-necessarily matching way, then coupling them together to arrive at a locally refined final mesh. This strategy allows for an easy local grid refinement in the vicinity of the structure while the mesh further away stays relatively coarse. In combination with the elaborate fine structure mesh (cf. Fig.\,\ref{Fig:Blocking}), this keeps the overall amount of elements relatively low while retaining well behaved element shapes. These blockwise meshes also easily allow to assign different polynomial degrees $p_i$ to the individual blocks.

\begin{figure}[h!]
    \hspace{-12cm}\textbf{Displacement:}\\
	\begin{center}
		\hspace{-2.5mm}\includegraphics[trim=0cm 0cm 0cm 0cm, clip, scale=0.35]{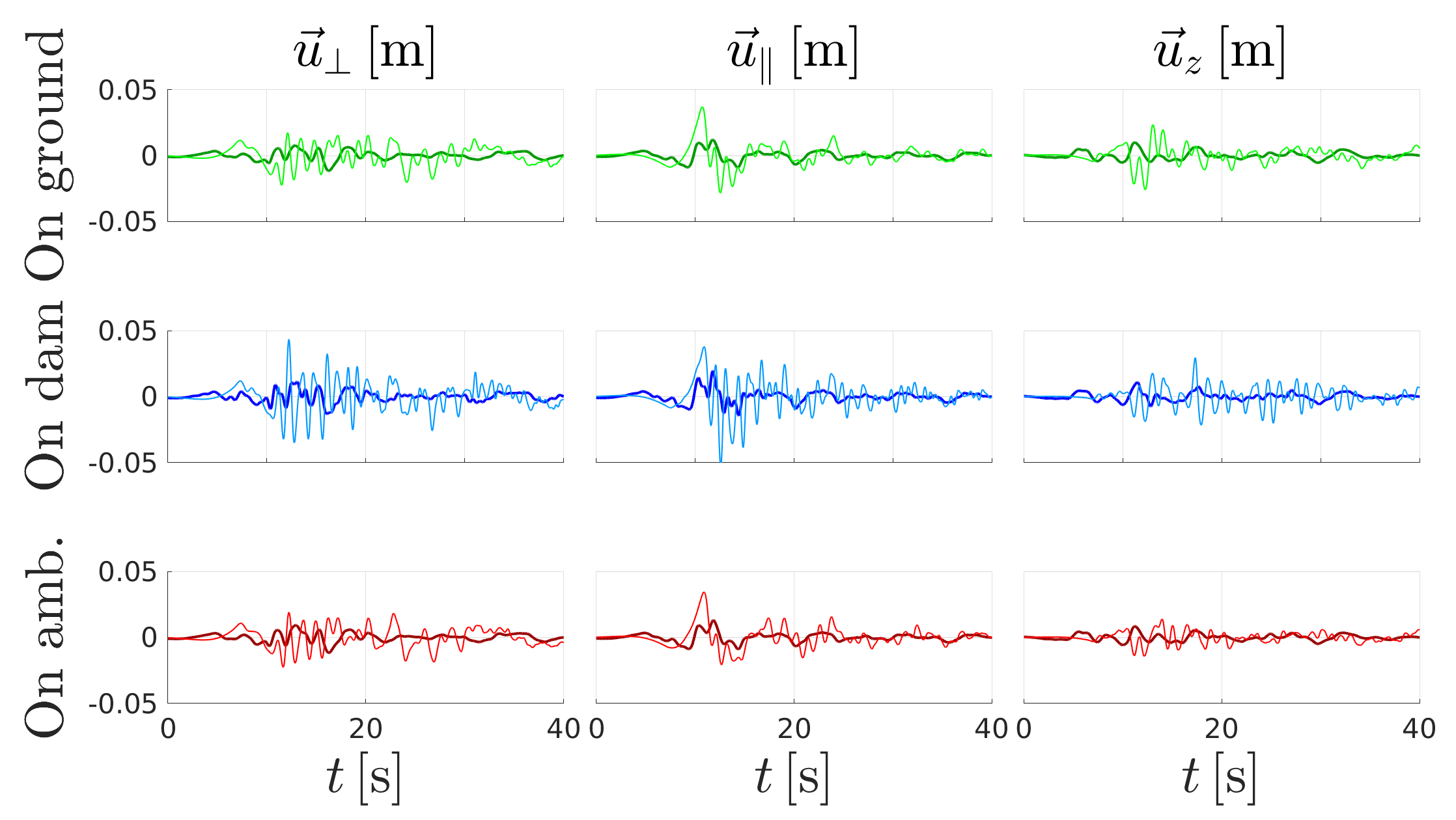}\\
    \end{center}
     \hspace{-13cm}\textbf{Velocity:}\\
    \begin{center}
		\includegraphics[trim=0cm 0cm 0cm 0cm, clip, scale=0.35]{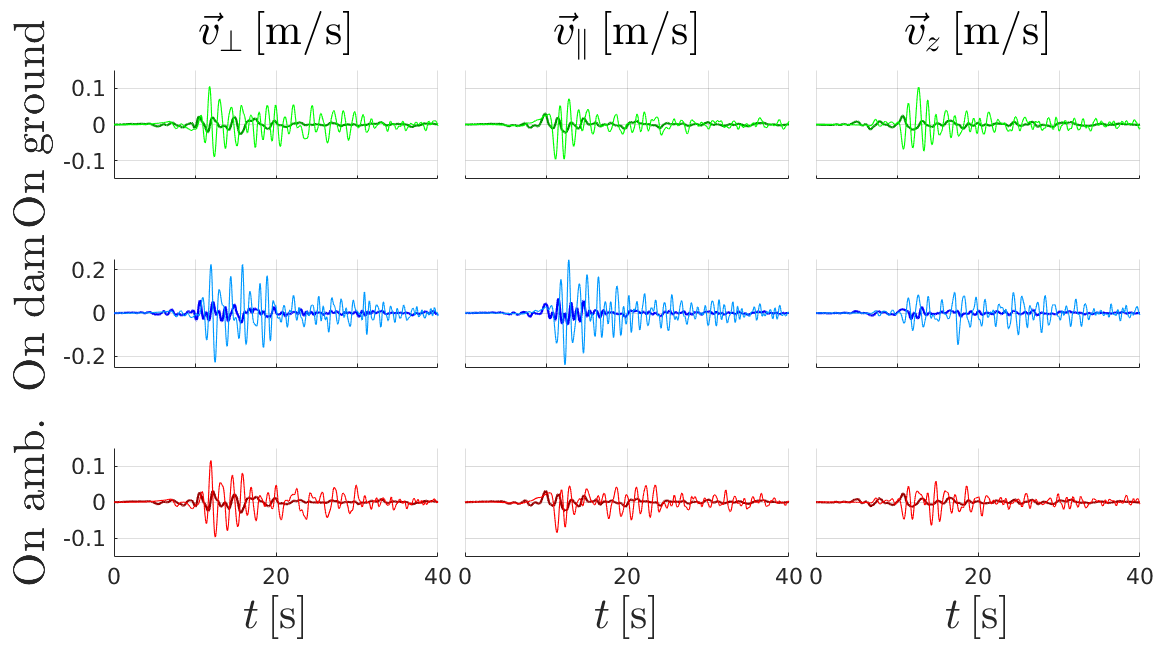}
	\end{center}
	\begin{center}
		\caption{\footnotesize Comparison of synthetical displacement \textbf{(top block)} and velocity \textbf{(bottom block)} seismograms between plain wave \textbf{(bold lines)} and source-to-site simulation \textbf{(thin lines)}. Each quantity is decomposed into dam-orthogonal, -parallel and vertical components and is compared at three physical locations (ground, dam-crest and ambutment), see Fig.\,\ref{Fig:UVA_Field} for locations. All data again filtered to $[0.1,1]\,\textup{Hz}$ with $2^{\textup{nd}}$ order Butterworth filter.\label{Fig:PlaneWave2vsS2S}}
	\end{center}
\end{figure}

\begin{figure}[h!]
	\begin{center}
	\hspace{-12cm}\textbf{Displacement:}\\
		\includegraphics[trim=0cm 0cm 0cm 0cm, clip, scale=0.15]{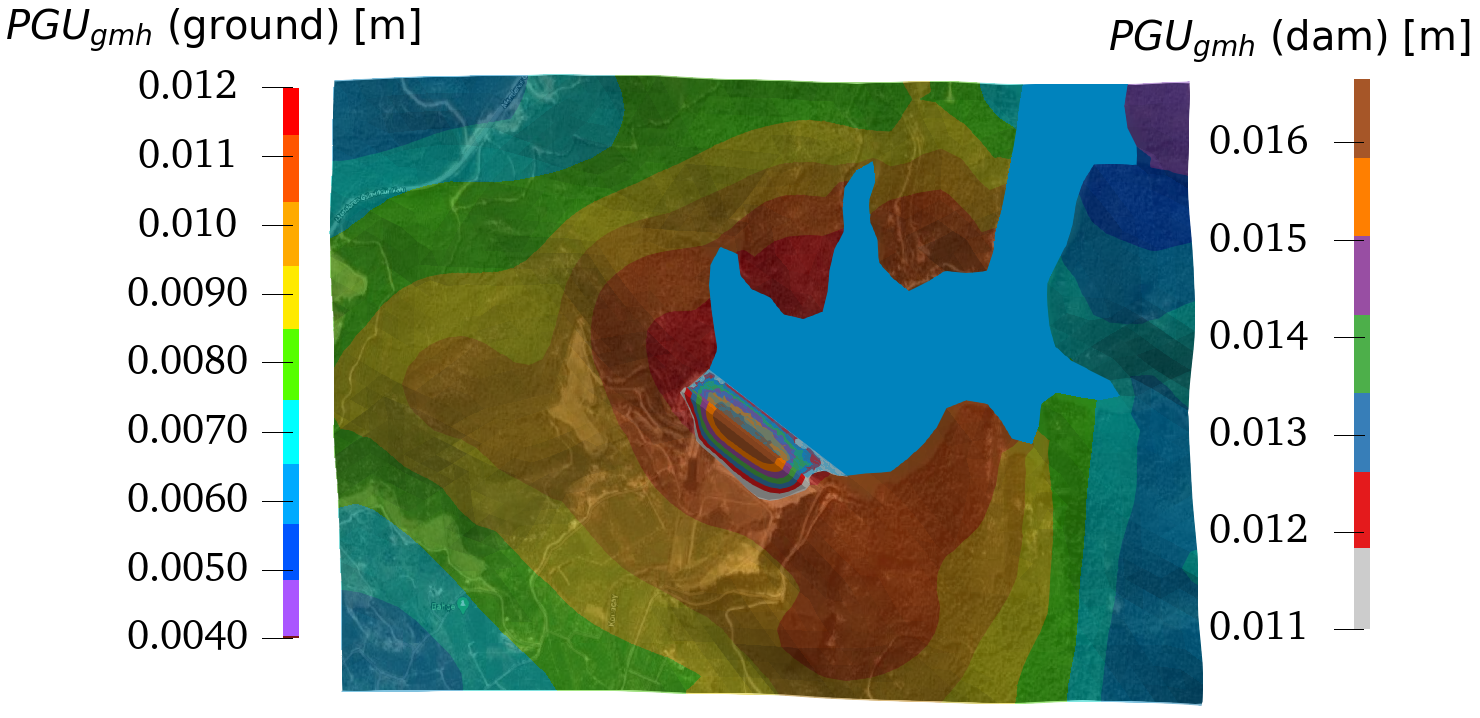}~~\includegraphics[trim=0cm 0cm 0cm 0cm, clip, scale=0.15]{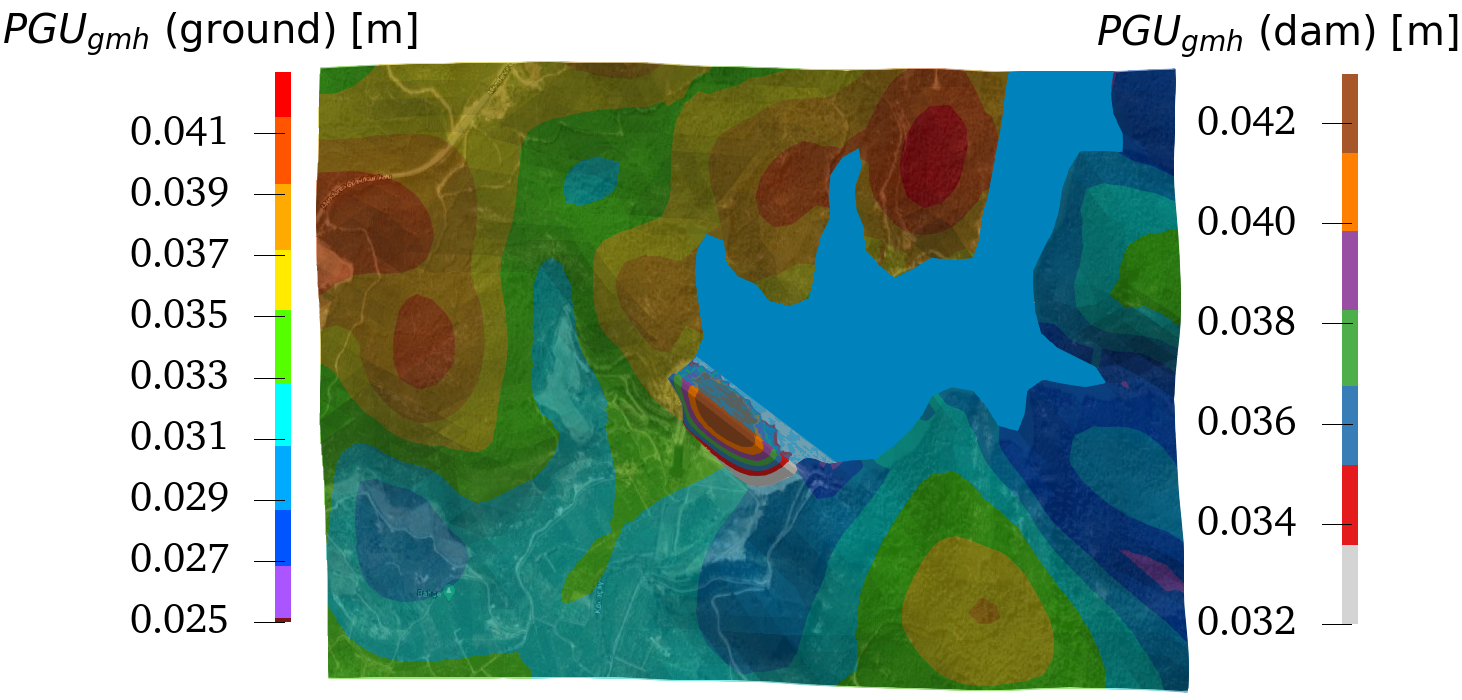}\\
    \hspace{-13cm}\textbf{Velocity:}\\	
		\includegraphics[trim=0cm 0cm 0cm 0cm, clip, scale=0.15]{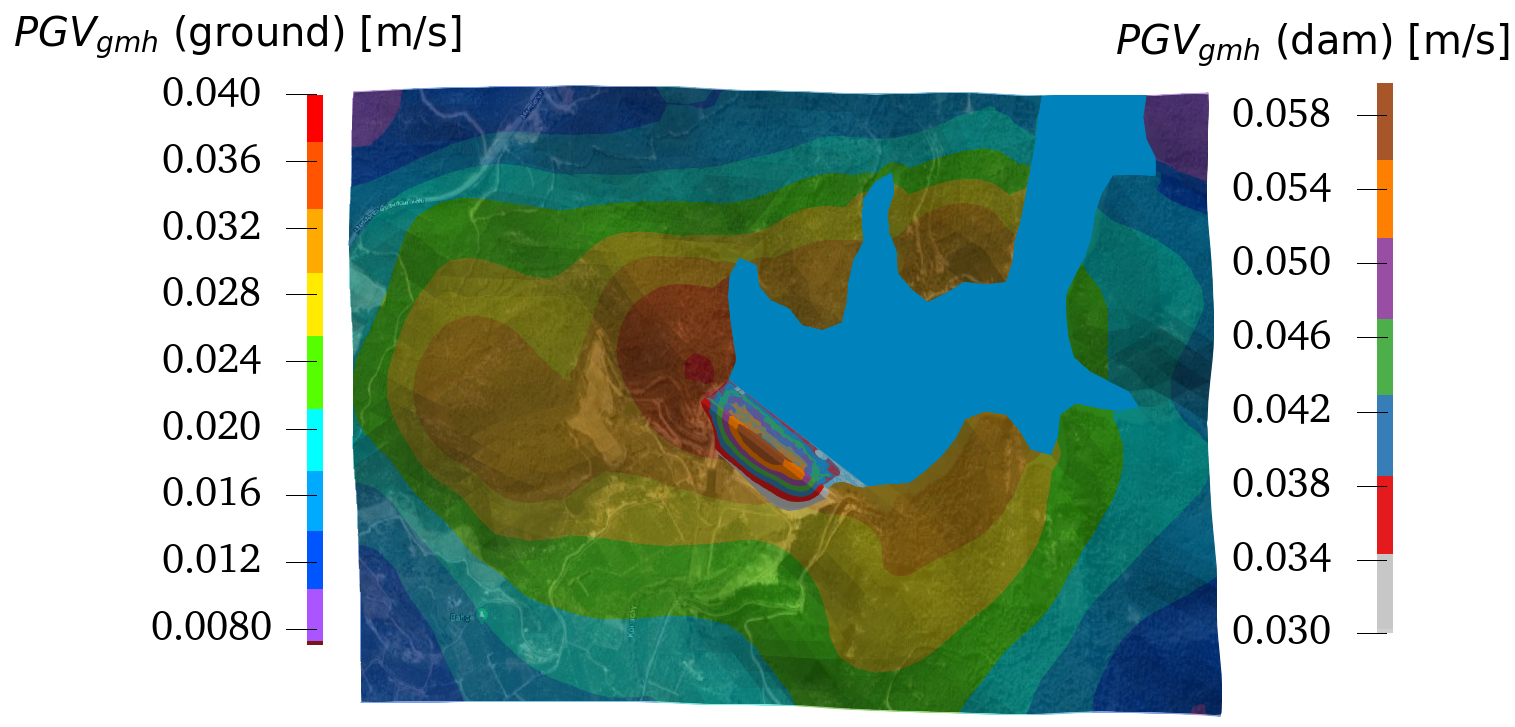}~~\includegraphics[trim=0cm 0cm 0cm 0cm, clip, scale=0.15]{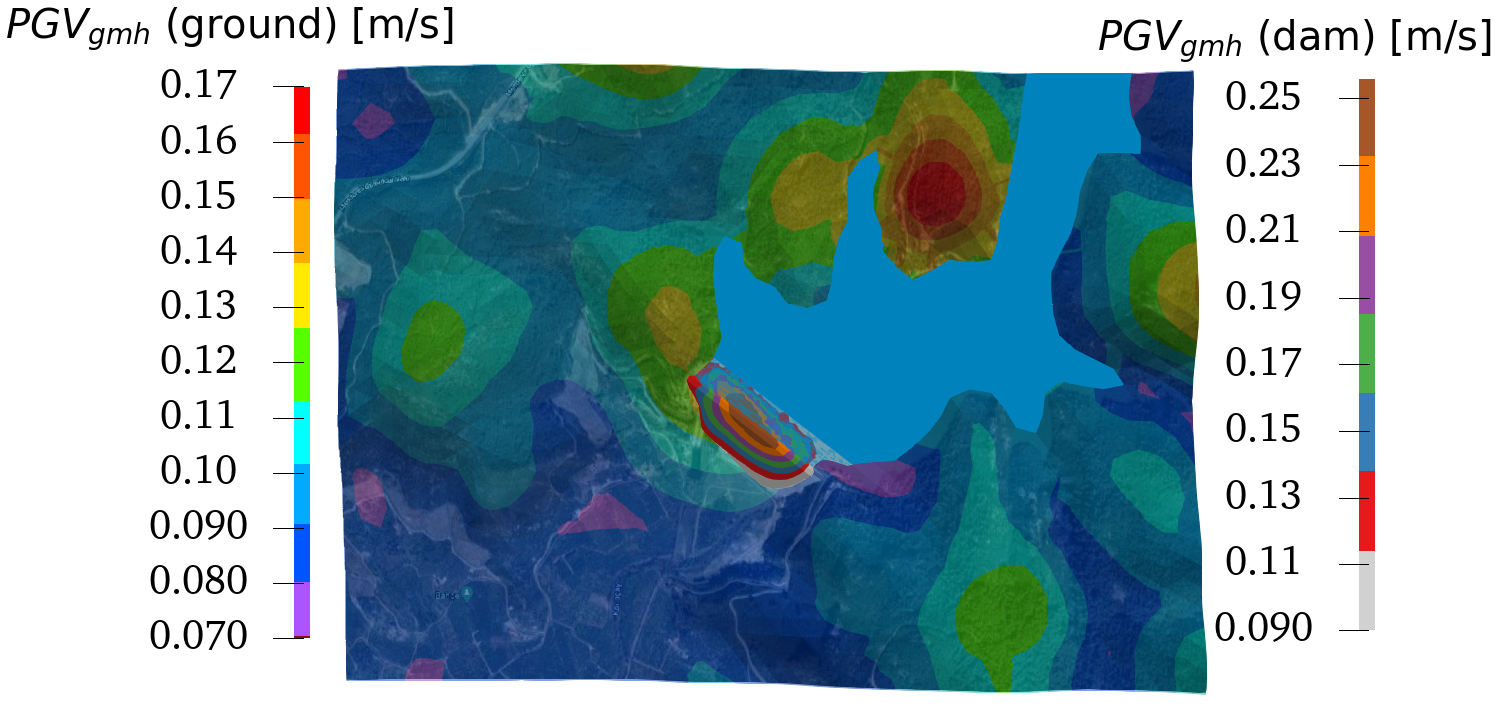}\\
		\vspace{-0mm}
		\hspace{9cm}\parbox{4.7cm}{\tiny Map data: Google, Imagery \textcopyright 2021 CNES/Airbus, Maxar Technologies, Map data \textcopyright 2021}
		\caption{\footnotesize Comparison of \textbf{(top row)} PGU and \textbf{(bottom row)} PGV maps between the \textbf{(left column)} plain wave simulation of Sec.\,\ref{subsec:PlaneWave1} and \textbf{(right column)} the full source-to-site simulation both in the vicinity region of the dam. \textit{(original map image overlayed with simulation colormap)}\label{Fig:PlaneWave2vsS2SPG}}
	\end{center}
\end{figure}

\begin{figure}[h!]
	\begin{center}
		\raisebox{4mm}{\includegraphics[trim=0cm 0cm 0cm 0cm, clip, scale=0.2]{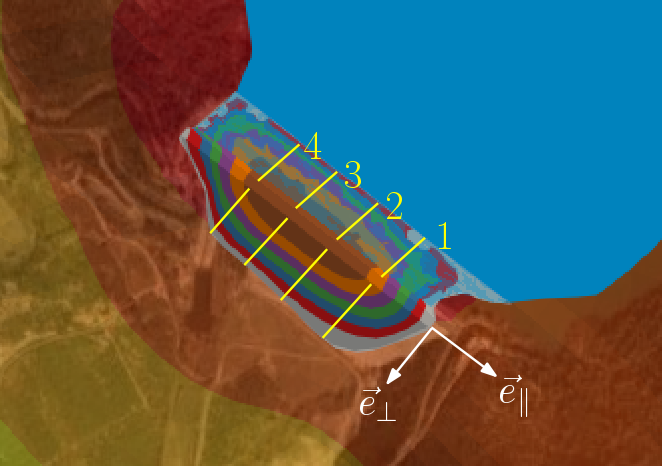}}\includegraphics[trim=0cm 0cm 0cm 0cm, clip, scale=0.3]{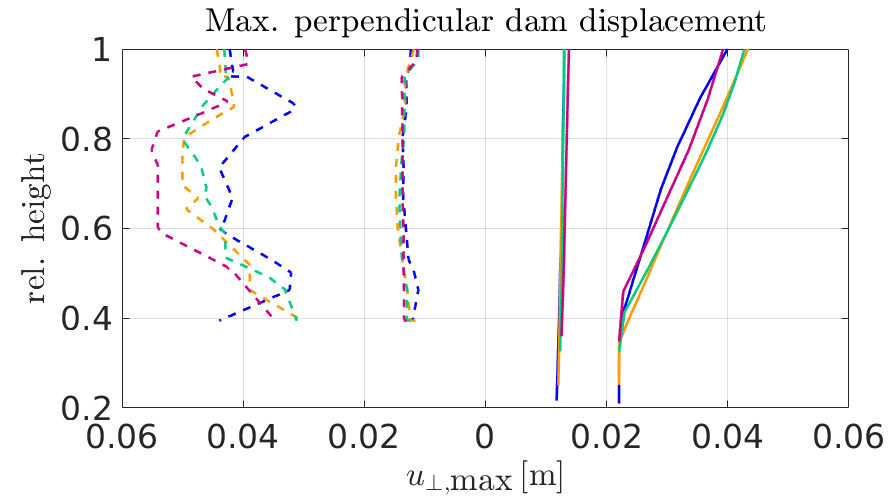}\\
		\vspace{-0mm}
		\hspace{-7cm}\parbox{4.7cm}{\tiny Map data: Google, Imagery \textcopyright 2021 CNES/Airbus, Maxar Technologies, Map data \textcopyright 2021}
		\caption{\footnotesize Evaluation of maximal orthogonal dam displacement $u_{\perp,\textup{max}}$ across four slices depicted in \textbf{(left)} image \textit{(original map image overlayed with simulation colormap)}. The image also shows the unit vectors in $\vec{e}_{\perp}$ in orthogonal and $\vec{e}_{\|}$ in parallel direction. Results show the \textbf{(middle strands)} plain wave simulation and \textbf{(outer strands)} source-to-site simulation with dashed lines for the reservoir side, continuous lines for the free side of the dam. ({\color{blue} blue}) slice 1, ({\color{orange} orange}) slice 2, ({\color{green} green}) slice 3, ({\color{magenta} red}) slice 4. \label{Fig:DamDisplacementSketch}}
	\end{center}
\end{figure}

As within the large scale simulation, also the source-to-site simulation adopts the layered materials, presented in Tab.\,\ref{Tab:MaterialParamLargeSim}, and the seismic fault rupture mechanism and slip distribution as previously described (see Fig.\,\ref{Fig:SourceData}). Both the layers and the location of the fault plane w.r.t. the dam structure are depicted in Fig. \ref{Fig:S2SDomain}.\\
The simulation was conducted using a polynomial degree of $p=2$ in the core- (water, dam, dam-vicinity) and the transition-blocks and $p=3$ in the outer layers and 920.000 timesteps of size $5\cdot 10^{-5}$. Fig.\,\ref{Fig:PlaneWave2vsS2S} shows the comparison between the synthetic seismograms computed on (i) ground, (ii) crest and (iii) abutment of the dam, as obtained with the source-to-site simulation (present paragraph) and with the plane-wave approach (previous simulation from Sec.\,\ref{subsec:PlaneWave1}). Note that the orthogonal, parallel and up-down synthetic time histories of displacement and velocity are considered; the directions are w.r.t to the orientation of the dam (see Fig.\,\ref{Fig:DamDisplacementSketch} for the orientations). Furthermore, Fig.\,\ref{Fig:PlaneWave2vsS2SPG} compares the respective peak ground maps of the two simulations: on the left-hand side the plane-wave model and, on the right-hand side, the source-to-site model. 


Finally Fig.\,\ref{Fig:DamDisplacementSketch} shows the maximum orthogonal displacement along the dam over different cross sections. The quantity $u_{\perp,\textup{max}}:=\sup_{t\in(0,T)}|\vec{u}_{\perp}|$, is presented along the four sections depicted in the nearby sketch.

\subsection{Discussion of results}
The results obtained so far deserve some comments, since at first glance the large difference in terms of displacement and velocity, experienced by the Tahtal{\i} dam adopting the plane-wave model and the source-to-site one, is evident. It is worth noting that the Tahtal{\i} dam and station AFAD \#3536 are located both at about 30\,km away from the hypocenter, therefore in the so-called near-field region, however the azimuthal difference is around 30$^{\circ}$. Radiation pattern \cite{Kotha2019} and rupture directivity \cite{Fukushima2003, Ripperger2008} effects might play a significant role in this region and therefore the large variability observed can be at least partially explained by these effects. Thanks to the validations accomplished (see Sec.\,\ref{subsec:Validation}) we are confident about the reliability of our simulations up to 1\,Hz, and therefore we consider the source-to-site simulation not only a state-of-the-art modelling approach but also the more reliable one, in terms of the excitation experienced by the Tahtal{\i} dam during the seismic event analysed in this study.\\
\indent Regarding the magnitude of the ground motion observed, it is important to mention that nor the plane-wave model neither the source-to-site simulation seems to be capable to produce shaking levels that might endanger the dam itself. In fact, according to the exhaustive literature examined, both based on numerical studies \cite{Anastasiadis2004, Bayraktar2009a, Bayraktar2009b,  Sunbul2017} or empirical observations \cite{seed1978performance, uscold1992observed, ussd2000observed, ussd2014observed, Yamaguchi2012, Park2018} the shaking level simulated seems to be incapable of producing significant damages to the infrastructure. This findings are coherent with the empirical observation as witnessed by the reconnaissance team \cite{Report_EERI,Report_METU}.\\
\indent Furthermore, as highlighted by \cite{Park2018}, there are more than 59,000 large dams worldwide, and more than three quarters employ earthfill and rockfill construction \cite{ICOLD2018}. Several large earthquakes were  recorded at embankment dams,  for example, during the 2008 Wenchuan earthquake ($M_w$ 7.9), the 156-meter-high Zipingpu concrete-faced rockfill dam (CFRD) was damaged partially without any collapse or freeboard deficiency. The dam, designed with peak ground acceleration of $0.26\,g$ at its foundation bedrock, recorded data exceeding $0.5\,g$ \cite{Kong2010, ZhangEtAl2015}. During the 2011 Tohoku earthquake ($M_w$ 9.0), the Aratozawa rockfill dam experienced a PGA of $1.04\,g$ at foundation rock and in spite of that the safety was not endangered \cite{JCOLD2014}. According to \cite{ussd2014observed} most modern embankment dams are capable of withstanding significant seismic shaking with no detrimental consequences in the past events. This leads to the conclusion that further analyses are necessary to predict more accurately which seismic event may involve a dangerous shaking for the structural safety of the building during earthquakes. The full source-to-site simulation provides an example of what our tool is capable of and must be properly exploited in the future.


\section{Conclusion and Outlook}
Starting from a general mathematical description of a coupled elasto-acoustic wave propagation problem, we have studied a realistic earthquake event for which we have analyzed the seismic response of a dam. The computational model comprises the actual topography around the dam, its reservoir lake as well as a simplified one dimensional crustal model. 
Regarding this last aspect, it is worth mentioning that, by taking into account a more accurate seismic tomography it will be possible to improve the computational model so far adopted. Eventual taking into account local soil heterogeneities in the proximity of the dam, as well as a more detailed characterisation of the material of the dam will also contribute to that. The generation and analysis of the latter will be the subject of future studies. Due to the comparably small ground motions recorded during the seismic event, we considered for the solid portion of the domain a relatively simple but rather realistic, visco-elastodynamic model; this latter should be enhanced to a plastic one, for an analysis focusing on individual features of the dam, especially in case of higher ground motions (i.e.: local events in the immediate proxmity of the dam). Having said that and precisely due to the model's simplicity, the conducted simulations turn out to be in good agreement with the recorded seismograms and are capable to produce reliable results in the frequency range of up to 1\,Hz at manageable computational costs. This resulted in a flexible and robust computational numerical model.\\
\indent The final, state-of-the-art, fully-coupled source-to-site simulation makes use of local and independent grid-refinements, treated with a discontinuous Galerkin approach, in order to accurately resolve the multiple length scales adopted in the model. With one single simulation the source-to-site approach allows to obtain numerical data at site, such as maximum displacement or peak velocities, that can be used for engineering purposes. The source-to-site model presented here could be used for a better assessment of the seismic risk associated to the dam and its nearby region by investigating the ground motion wave field generated by (i) different earthquake realizations along the Kaystrios fault and/or (ii) different seismogenic faults. Due to its generality the model can be also easily employed and further adapted for the seismic risk assessment in other active regions and for different structures with additional uncertainties being considered.


\section*{Acknowledgements}
M. Muhr and B. Wohlmuth acknowledge the financial support provided by the Deutsche Forschungsgemeinschaft under the grant number WO 671/11-1. I. Mazzieri is member of the INdAM Research group GNCS and this work is partially funded by INdAM-GNCS.

\footnotesize
\bibliographystyle{abbrv} 
\bibliography{references}

\end{document}